\documentclass[10pt]{amsart}
  \usepackage{amsmath, amsfonts, amssymb, amsthm, latexsym, mathrsfs, eufrak}
  \usepackage[all]{xy}

\def\A{\mathcal A}
\def\DD{\mathcal D}
\def\F{\mathcal F}
\def\H{\mathcal H}
\def\I{\mathcal I}
\def\L{\mathcal L}
\def\O{\mathcal O}

\def\V{\mathcal V}
\def\P{\mathbb P}
\def\C{\mathbb C}

\def\mm{\mathfrak{m}}
\def\pp{\mathfrak{p}}

\def\ft{\tilde \varphi}
\def\htil{\tilde h}
\def\kt{\tilde k}
\def\ltil{\tilde l}
\def\Phit{\tilde \Phi}

\def\st{\tilde \sigma}
\def\Zt{\tilde Z}
\def\ab{\bar a}
\def\bb{\bar b}
\def\jb{\bar{\text{\emph{\j}}}}
\def\lb{\bar l}
\def\qb{\bar q}
\def\sb{\bar s}

\def\vb{\bar v}

\def\xb{\bar{\xi}}
\def\zb{\bar z}
\def\Xb{\overline{X}}

\def\fo{\text{f.o.}}
\def\lra{\longrightarrow}
\def\ra{\rightarrow}
\def\implies{\Longrightarrow}
\def\isom{\simeq}
\def\ge{\geqslant}
\def\le{\leqslant}
\def\tensor{\otimes}
\def\de{\partial}
\def\deb{\bar{\partial}}
\def\i{\sqrt{-1}}

\let\a\alpha
\let\b\beta
\let\g\gamma
\let\d\delta
\let\e\epsilon
\let\ep\varepsilon
\let\z\zeta
\let\x\xi
\let\th\theta
\let\del\vartheta
\let\k\kappa
\let\l\lambda
\let\m\mu
\let\n\nu
\let\p\pi
\let\r\rho
\let\s\sigma
\let\t\tau
\let\f\varphi
\let\o\omega
\let\G\Gamma
\let\D\Delta
\let\Th\Theta

\let\Si\Sigma

\begin{document}
\title{{On Vector Bundles of Finite Order}}
\author{Mario Maican}
\address{Department of Mathematics, University of California,
         Riverside, CA 92521}
\email{mmaican@math.ucr.edu}

\maketitle

\tableofcontents


\section*{Introduction}

Bundles of finite order were first introduced in \cite{griffiths}
and were studied systematically in \cite{griffiths-cornalba}.
Griffiths' ``working hypothesis'' is that on an affine variety
(smooth, over $\C$ and satisfying some extra conditions, cf. \S 1)
every holomorphic vector bundle has a finite order structure.
This he calls the ``Oka principle with growth conditions''
by analogy with the Oka-Grauert principle which states that
on a Stein manifold every topological vector bundle has a unique
holomorphic structure. 

Griffiths and Cornalba prove the Oka
principle for line bundles only. See (4.4) for the precise
statement. For bundles of rank greater than 1 their theory is
somewhat insufficient. To understand why we review here the four
definitions they gave for what should mean that a holomorphic
vector bundle $E$ on $X$ has finite order: 

\begin{enumerate}
\item[(I)] $E$ has a Hermitian metric whose holomorphic bisectional
curvature has polynomial growth, see (4.1);
\item[(II)] there is a holomorphic map $f$ of finite order
from $X$ to some Grassmannian such that $E\isom f^* \mathcal{U}$, where
$\mathcal{U}$ is the universal bundle on the Grassmannian;
\item[(III)] $E$ has Schubert cycles of finite order;
\item[(IV)] $E$ has transition matrices of finite order relative to
a suitable covering of $X$ with punctured polycylinders.
\end{enumerate}

\noindent
We refer to section 2 for the notions of holomorphic map of
finite order, resp. analytic subset of finite order.
The polycylinders from (IV) are as in \S 1. In the case of a line
bundle $L$ the above conditions take a simpler form:

\begin{enumerate}
\item[(I)] $L$ has a Hermitian metric whose first Chern form has 
polynomial growth;
\item[(II)] there is a holomorphic map of finite order $f:X\lra \P^N$
such that $L\isom f^* \O_{\P^N}(1)$;
\item[(III)] $L=[Z]$ with $Z\subset X$ an analytic divisor of finite order;
\item[(IV)] $L$ has transition functions of finite order relative
to a suitable covering of $X$. 
\end{enumerate}

As shown in \cite{griffiths-cornalba} the four conditions from above are
equivalent in the case of a line bundle.
The implication ``(I)$\implies$(II)''
is true for any rank and follows from the vanishing theorem (4.6).

The implication ``(II)$\implies$(III)'' fails in general
because of the so-called
``unsolvability of the Bez\^out problem in codimension greater than 1''.
In a form tailored to our purposes, the Bez\^out problem asks whether
intersection of sets of finite order is of finite order, too.
It has negative answer because of the counter-example from
\cite{cornalba-shiffman}. Shiffman and Cornalba give two analytic sets
of order zero in $\C^2$ whose intersection has infinite order.

The implication ``(III)$\implies$(IV)'' cannot be carried out
if the rank is greater than 1 because of
technical reasons: when $k\ge 2$ one does not know whether an
analytic set of codimension $k$ in $\C^n$ can be given as the
zero-set of $k$ functions of finite order. The case $k=1$ has positive
answer, cf. \cite{skoda}. In the case $k>1$ Skoda showed that one
can find $n+1$ defining functions of finite order, cf. \cite{skoda_2},
but this is of little help.

In view of these difficulties we propose the following approach:
one should give the conditions for finite order growth not in terms
of $E$ but in terms of the hyperplane bundle $\O_{\P(E^*)}(1)$.
This approach is inspired from classical complex geometry:
one knows that a holomorphic vector bundle $E$ over a compact
complex manifold is ample (in the sense that the zero-section
of its dual can be collapsed to a point, or in the sense that Cartan's
Theorem A or Theorem B hold) if and only if $\O_{\P(E^*)}(1)$
is ample. Denoting $\L= \O_{\P(E)}(1)$ we propose the following definition
for what should mean that the bundle $E^*$ have finite order:

\begin{enumerate}
\item[(I')] $E$ has a Finsler metric of finite order, cf. (5.5).
\end{enumerate}

\noindent
One is forced to think about Finsler metrics because giving
a Finsler metric on $E$ is, roughly speaking, equivalent to giving
a Hermitian metric on $\L$. For Finsler metrics one has the notion
of holomorphic bisectional curvature, cf. (3.21), and horizontal holomorphic
bisectional curvature, cf. (3.22). We define finite order Finsler metrics by
imposing estimates on the curvature in a very similar manner
to the Hermitian case.

Once we are given a Finsler metric on $E$ there is a natural way of doing
Nevanlinna theory on $\P(E)$ which we explain in section 7.
Thus we can formulate a definition analogous to (II):

\begin{enumerate}
\item[(II')] there is a holomorphic map of finite order to some
projective space $f:\P(E)\lra \P^N$ such that $f^* \O_{\P^N}(1)\isom \L$.
\end{enumerate}

The crucial step in the proof of ``(I)$\implies$(II)'' is a vanishing
theorem for the sheaf of sections of finite order of $E$.
At (5.12) we prove the corresponding statement in a Finsler context:

\noindent \\
{\bf Theorem 1:} \emph{Let $X$ be a special affine variety
and $E$ a holomorphic
vector bundle on $X$ equipped with a Finsler metric of finite order.
Then, for any $q\ge 1$, the $q^{\text{th}}$ cohomology of $\O_{\fo}(E^*)$
vanishes.}

\noindent \\
As a consequence we get at (7.10) the implication ``(I')$\implies$(II')'':

\noindent \\
{\bf Theorem 2:} \emph{Let $X$ be a special affine variety of dimension $n$.
Let $E$ be a holomorphic vector bundle of rank $r+1$ on $X$ equipped
with a Finsler metric of finite order. Then, for any integer
$N > 2(n+r)$, there is a holomorphic immersion $f:\P(E) \lra \P^N$
of finite order satisfying $f^* \O_{\P^N}(1) \isom \L$.}\\

In \S 1 we introduce the spaces on which we are working.
They come equipped with complete K\"ahler metrics of finite volume
and bounded Ricci curvature. In \S 2 we present just the amount of
Nevanlinna theory that we need. In \S 3 we define Finsler metrics
and, following \cite{cao-wong}, we compute the curvature of $\L$.
In \S 4 we recall the main results of Griffiths and Cornalba on
Hermitian metrics of finite order. We begin to make the transition
to Finsler metrics by giving an equivalent definition for sections
of finite order, cf. (4.12). \S 5 and \S 7 contain our two
main results as stated above. They should be understood as a
very small step in attempting to elucidate the Oka principle with
growth conditions.
In \S 6 we apply the classical theory
of embeddings of Stein manifolds to show that the map from
theorem 2 can be made to be an immersion. Notice that we do not
claim that the image of $f$ in $\P^N$ is closed and, in fact,
it is not. \\ \\

%
%

\section{Preliminaries}

We begin by introducing the spaces on which we will be working.
Let $\overline{X}$ be a complex projective manifold of dimension $n$.
Let $D_1,\ldots, D_\n$ be effective smooth ample divisors on $\Xb$.
We put $D= D_1+\ldots +D_\n$. We assume that $D$ has simple normal
crossings. This means that around each point $x\in \Xb$ there is a
coordinate chart $(z_1,\ldots, z_n)$ in which $D =\{ z_1 \cdots
z_k =0\}$ for some $0\le k\le n$. We put $X= \Xb \setminus D$
and call $X$ a \emph{special affine variety}. Note that, by Hironaka's
resolution of singularities, a smooth affine variety always admits a
compactification such that the divisor at infinity have simple normal
crossings. Thus the word ``special'' refers to our requirement that
each $D_i$ be ample.

Let $\D = \{ z\in \C,\ |z|\le  1\}$ be the unit disc and $\D^* =
\{ z\in \C,\ 0< |z| \le 1 \}$ be the punctured unit disc. A
\emph{k-fold punctured polycylinder} is of the form $P^* = (\D^*)^k
\times \D^{n-k},\ \ 0\le k \le n$. The compact manifold $\Xb$ can be
covered with finitely many polycylinders of the form $P= \D^n$.
They can be chosen in such a manner that the intersections
$P^* =P\cap X$ be k-fold punctured polycylinders. If $k=0$ then we
have a polycylinder that is entirely included in $X$. If $k>0$
we call $P^*$ a \emph{neighbourhood at infinity}.

On the punctured unit disc we have the Poincar\'e metric
\begin{align*}
\o_{\D^*}= \frac{\i}{\p}\cdot \frac{dz \wedge d \bar{z}}{|z|^2
(\text{log} \frac{|z|^2}{c})^2} = - dd^c \text{log}(\text{log} 
\frac{|z|^2}{c})^2 ,\ \ c>1,
\end{align*}
which has constant negative Gauss curvature, is complete and
has finite volume. By the Poincar\'e metric $\o_{P^*}$ on a punctured
polycylinder we mean the product of Poincar\'e metrics on the punctured
components and Euclidean metrics on the non-punctured components.
For a positive (1,1)-form $\f$ on $X$ we write
\begin{align*}
\tag{1.1}
\f \sim \o_{P^*}
\end{align*}
if on each polycylinder at infinity $\f$ is equivalent to the
Poincar\'e metric.

Using the assumption that $D_i$ are ample Griffiths and Cornalba
construct an exhaustive function $\t$ on $X$ with the following
properties:

\begin{enumerate}
\item[(i)] $\t$ is strictly plurisubharmonic;
\item[(ii)] the Levi-form $\f = dd^c \t$ induces a complete
metric on $X$;
\item[(iii)] the Ricci curvature of this metric is bounded:
\begin{align*}
\tag{1.2}
|\text{Ric}(\f)| = O(\f);
\end{align*}
\item[(iv)] $\f$ has finite volume.
\end{enumerate}
The contruction of $\t$ goes as follows:
Consider the line bundles $[D_i]$ on $\Xb$ associated to $D_i$. They
are ample, so on each of them we can choose a Hermitian metric with
positive Chern form. We choose global sections $\s_i$ of $[D_i]$ whose
zero-sets coincide with $D_i$ and which have norm less than 1 at every
point. We put $\s= \s_1 \tensor \ldots \tensor \s_\n$, which is a
section of $[D]= [D_1]\tensor \ldots \tensor [D_\n]$, and define
\begin{align*}
\tag{1.3}
\t=c\ \text{log} \frac{1}{|\s |^2} - \
\text{log} \{ (\text{log}|\s_1 |^2)^2 \cdots 
(\text{log}|\s_\n |^2)^2 \}.
\end{align*}
Notice that $dd^c$log$\frac{1}{|\s |^2} = c_1 ([D])$ is a positive
(1,1)-form on $\Xb$, while \\
$dd^c$ log(log$ |\s_i |^2)^2$ is the
Poincar\'e metric in the vicinity of $D_i$. Thus, choosing $c$ large
enough, we can make sure that the Levi-form $\f =dd^c \t$ be
positive, i.e. (i) holds. In addition $\f \sim \o_{P^*}$ which
tells us that (ii) and (iv) hold. The remaining property (iii)
was proven in \cite{griffiths-cornalba} and amounts to the
fact mentioned above that $\o_{\D^*}$ has bounded curvature.
We shall also make use of the exhaustive function $\r = e^{\t/2}$.

From now on $X$ will be a special affine variety. We fix a compactification
$\Xb$ and an exhaustive function $\t$ as above. We assume that $X$ is
embedded in some $\C^m$ and $\Xb$ is the closure of $X$ in $\P^m$. Such
an embedding exists because the line bundle $[D]$ associated to the divisor
at infinity is ample. In addition, we may assume that $[D]$ is the
restriction of $\O_{\P^m}(1)$ to $X$. We choose affine coordinates
$z=(z_1,\ldots, z_m)$ on $\C^m$ and homogeneous coordinates
$(z_0;\ldots; z_m)$ on $\P^m$. As global section of $[D]$ in (1.3) we can
take $\s = z_0$. We have $|\s |^2 = (1+||z||^2)^{-1}$ and $\r \sim 
(1+||z||^2)^{c/2}$ outside a compact subset of $X$. Hence
\begin{align*}
\tag{1.4}
\r \sim ||z||^c
\end{align*}
outside a compact subset of $X$.

Next we introduce a semipositive (1,1)-form $\psi$ on $X$ which is suitable
for doing Nevanlinna theory. We fix a projection $\pp : X\lra \C^n $ onto
an n-dimensional subspace of $\C^m$ such that $\pp$ is a branched covering.
We put
\begin{align*}
\t'= \text{log} ||\pp (x)||^2,\ \ \r'=e^{\t'/2},\ \ \psi = dd^c \t'.
\end{align*}
We have $\psi \ge 0,\ \psi^{n-1}\neq 0$ and $\psi^n =0$.
Thus $\psi$ is not positive definite. We notice that outside a compact
subset we have
\begin{align*}
\tag{1.5}
\t \sim \t'.
\end{align*}

We finish this section with several estimates from \cite{maican}
that we will need later:
\begin{align*}
\tag{1.6}
\psi \le \r^c \f,
\end{align*}
\begin{align*}
\tag{1.7}
d\t' \wedge d^c \t' \le \r^c \f,
\end{align*}
\begin{align*}
\tag{1.8}
d\t \wedge d^c \t \le \r^c \f
\end{align*}
for some positive constant $c$. \\ \\

%
%

\section{Nevanlinna Theory on Special Affine Varieties}

In this section we will recall several standard notions
from Nevanlinna theory: the characteristic function measuring
the growth of a holomorphic map from a special affine variety
$X$ to projective space, the counting function measuring the
growth of an analytic subset of $X$. We will then formulate
basic relationships between these functions. As general reference
we point out the elegantly written \cite{shabat}.

Relative to the exhaustive function $\r' = e^{\t'/2}$ introduced
in \S 1 we put
\begin{align*}
X[r]= \{ x\in X,\ \r'(x) \le r\},\quad X<r>=\{x\in X,\ \r'(x)= r\}.
\end{align*}
By Sard's theorem the sets $X<r>$ are smooth for all $r$ outside a set of
measure zero. In the sequel, each time we integrate over the set $X<r>$,
it will tacitly be assumed that the latter is smooth.

Given a holomorphic map $f:X\lra \P^N$ we define its \emph{characteristic
function}
\begin{align*}
T_f (r,s) = \int_s^r \frac{dt}{t} \int_{X[t]} f^* \o \wedge \psi^{n-1}
\end{align*}
and its \emph{higher characteristic functions}
\begin{align*}
T_f^{(k)} (r,s) = \int_s^r \frac{dt}{t} \int_{X[t]} f^* \o^k \wedge \psi^{n-k},
\end{align*}
where $r>s>0$ are real numbers, $1\le k \le n$ and $\o$ is the
Fubini-Study form on $\P^N$. Given an analytic subset $Z\subset X$ of
pure dimension $k$ we define its \emph{counting function}
\begin{align*}
N_Z (r,s) = \int_s^r \frac{dt}{t} \int_{Z[t]} \psi^k,
\end{align*}
where $Z[t]=Z\cap X[t]$. In this definition and hereafter $Z$ does not
have to be reduced, in other words its components may have multiplicities.
Note that if the image of $Z$ under the projection
$\pp :X\lra \C^n$ does not contain the origin then, the quantity $N_Z(r)=
N_Z (r,0)$ is well-defined. We call it ``counting function'' because
in the case $k=0$, i.e. when $Z$ is a sequence of points, $N_Z (r)$
equals the logarithmic average of the number of points of $Z$
inside the ball $X[r]$.
Given a global section $\s$ of $\O_{\P^N}(1)$
we define the \emph{proximity function} of $f$ to the zero-set of $\s$ by
\begin{align*}
m_\s (r) =\int_{X<r>} \text{log} \frac{1}{|\s \circ f|} d^c \t' \wedge
\psi^{n-1}.
\end{align*}
Here $|\cdot |$ is the canonical norm of $\O_{\P^N}(1)$ given at a point
$[v]$ by $|\s |_{[v]} =\frac{|<\s, v>|}{|v|}$. Notice that $d^c \t' \wedge
\psi^{n-1}$ is a volume form on $X<r>$. We call $m_\s (r)$ ``proximity
function'' because it takes large values whenever the image of $X<r>$
under $f$ is close to the zero-set of $\s$. By choosing $\s$ to have norm
less than 1 at all points we can arrange that $m_\s (r)$ be non-negative.\\

\noindent
{\bf (2.1) First Main Theorem:} 
\emph{Let $f:X\lra \P^N$ be a holomorphic map.
Let $\s$ be a global section of $\O_{\P^N}(1)$ and
$Z= f^* \{ \s =0\}$.
Assume that the image of $f$ is not contained in $\{ \s = 0 \}$.
Then, for $r>s>0$, we have}
\begin{align*}
N_Z(r,s)+m_\s(r)-m_\s(s)= T_f(r,s).
\end{align*}

\noindent
{\bf (2.2) Corollary} (Nevanlinna's Inequality):
\emph{Fix $s>0$. Then, under the
hypotheses of the previous theorem, we have for $r>s$}
\begin{align*}
N_Z(r,s)\le T_f(r,s) +O(1).
\end{align*}
\emph{Here $O(1)$ is a constant that may depend on $s$.}

\noindent \\
Thus the growth of the preimage of a hyperplane is bounded by the growth
of $f$. The corresponding statement for planes of codimension greater than 1 is
false: Cornalba and Shiffman give in \cite{cornalba-shiffman} an example of
a holomorphic map of order zero from $\C^2$ to $\C^2$ for which the
preimage of the origin has infinite order. This falls into the circle
of ideas known as ``Bez\^out problem''. More about the Bez\^out problem
we said in the introduction. We now state a converse to (2.2) in an
``average'' sense. To explain it we introduce the Grassmannian
$G(N,k)$ of planes of codimension $k$ in $\P^N$.
There is a unique measure $\m$ on $G(N,k)$ which is invariant under
the action of the unitary group $U(n+1)$ and which is so normalized
that the measure of the total space be 1.

\noindent \\
{\bf (2.3) Crofton Formula:}
\emph{Let $f:X\lra \P^N$ be a holomorphic map
which is non-degenerate in the sense that the preimage of a plane $P$ of
codimension $k$ in $\P^N$ is an analytic subset of codimension $k$ in $X$.
Then, for $r>s>0$, we have}
\begin{align*}
T_f^{(k)} (r,s)= \int_{P\in G(N,k)} N_{f^* P} (r,s) d\m (P).
\end{align*}

\noindent \\
{\bf (2.4) Definition:} Let $f:X\lra \P^N$ be a holomorphic map. 
We say that $f$ \emph{has finite order} if there is $\l \ge 0$ such
that
\begin{align*}
T_f (r,s)= O(r^\l).
\end{align*}
Likewise, we say that an analytic subset $Z\subset X$ of pure dimension
$k$ \emph{has finite order} if there is $\l \ge 0$ such that
\begin{align*}
N_Z (r,s)= O(r^\l).
\end{align*}
By an abuse of language we will say that $\l$ \emph{is the order of} $f$
or of $Z$ if the above estimates hold. Technically, we would have to
say that ``$Z$ has order at most $\l$'' but we want to avoid obstinate
repetitions.

\noindent \\
{\bf (2.5) Remark:} Assume that $f$ has finite order and is linearly
non-degenerate. Then (2.2) tells us that for any hyperplane $H\subset
\P^N$ the preimage $f^* H$ has finite order. Conversely, assume that the
preimages $f^* H$ have finite order in a uniform fashion, i.e. there
are $r_o, \k, \l \ge 0$ such that for $r\ge r_o$ and all $H$ we have
$N_{f^* H} (r,s)\le \k r^\l$. Then, by (2.3), also $f$ has finite order.\\

The order of growth of an analytic subset $Z \subset X$ depends
on the embedding $X\subset \C^m$ and the projection $\pp :X\lra \C^n$
that we fixed in \S 1. However, the notion of $Z$ having finite order
is intrinsic. We see this as follows: assume that $Z$ has pure dimension
$k$ and consider the counting function of $Z$ computed using $\f$
instead of $\psi$:
\begin{align*}
\hat{N}_Z (r,s) = \int_s^r \frac{dt}{t} \int_{Z_{\r \le t}} \f^k.
\end{align*}
Using (1.5), (1.7) and the technique from \cite{griffiths-carlson}
we can show that $N_Z (r,s) = O(r^\l)$ if and only if there is $\m$
depending only on $\l$ such that $\hat{N}_Z (r,s)= O(r^\m)$.
See \cite{maican} for the details.
Now, if $\t_1$ and $\t_2$ are constructed starting from two different
embeddings of $X$ then it is an easy matter to see that $\t_1 \sim \t_2$
and $dd^c \t_1 \sim dd^c \t_2$. Thus $\hat{N}_1$ and $\hat{N}_2$ have
polynomial growth at the same time.

By (2.3) the same considerations apply to holomorphic maps
$f:X\lra \P^N$. The order of growth of $f$ is not intrinsic however,
the notion of $f$ having finite order does not depend on any of the
choices made.

Finally, we would like to mention the sheaf of holomorphic
functions of finite order. First note that a holomorphic function
$f:X\lra \C$ can be regarded as a map to $\P^1$ by sending $x$ to
$(1;f(x))$. Its characteristic function takes the form
\begin{align*}
T_f (r,s)= \int_s^r \frac{dt}{t} \int_{X[t]} dd^c \text{log}(1+|f|^2)
                \wedge \psi^{n-1}.
\end{align*}
Let us introduce one last growth function, the \emph{maximum modulus
function}
\begin{align*}
M_f (r)= \text{log max}\{ |f(x)|,\ \r(x)=r\}.
\end{align*}
Using the fact that log$(1+|f|^2)$ is plurisubharmonic one
can show that $T_f (r,s)$ has roughly the same growth as $M_f (r)$.
This point of view allows us to localize the notion of finite order
function:

\noindent \\
{\bf (2.6) Definition:} We define the \emph{sheaf $\O_\l$ of germs of
holomorphic functions of order $\l$ on $X$} as follows: $\O_\l$ is a
sheaf on $\Xb$; for each open set $U\subset \Xb$ the space of sections
$\O_\l (U)$ consists of those holomorphic functions on $U\cap X$
with the property that around each point at infinity $x\in (\Xb \setminus
X)\cap U$ there is a neighbourhood $W\subset U$ and a constant $\k\ge 0$
such that the estimate $|f(w)|\le$ exp$(\k r^\l)$ with $w\in W,\
r=\r (w)$ holds.

Similarly we define the \emph{sheaf $\O_{\fo}$ of germs of holomorphic
functions of finite order on $X$} by the requirement that the estimate
$|f(w)|\le$ exp$(\k r^\l)$ hold with $\k, \l$ depending on the point
at infinity $x$.

It is easily seen that $\O_\l$ and $\O_{\fo}$ are sheaves on $\Xb$.
In fact they are $\O_{\Xb}$-modules. Their restrictions to $X$ coincide
with the sheaf $\O_X$ of germs of holomorphic functions on $X$.
By the compactness of the divisor at infinity it is also clear that
the spaces of global sections $\O_\l (\Xb)$ and $\O_{\fo} (\Xb)$
coincide with the space of holomorphic functions of order $\l$ on
$X$, respectively the space of holomorphic functions of finite order
on $X$.

These sheaves were studied in \cite{griffiths-cornalba} and
\cite{wong-mulflur}. Griffiths and Cornalba showed that $\O_\l$
and $\O_{\fo}$ are acyclic, i.e. their cohomology groups, apart
from $H^0$, vanish. Wong et al. proved that $\O_{\fo}$ is flat over
$\O_{\Xb}$ and as a corollary obtained the acyclicity of $\O_{\fo}$.\\

%
%

\section{Finsler metrics and the geometry of $\P(E)$}

Let $X$ be a special affine variety of dimension $n$.
Let $E$ be a holomorphic vector bundle of rank $r+1$ on
$X$. We begin by recalling some facts about Hermitian
metrics on $E$. We will then introduce Finsler metrics
and see how the notion of holomorphic bisectional curvature
generalizes from the Hermitian case to the Finsler case. Most of
our calculations are taken from \cite{cao-wong}. We also
refer to \cite{abate}.

We will denote by $\O_X, \A_X, \A_X^k, \A_X^{p,q}$ the sheaves of
holomorphic functions, of smooth $\C$-valued functions, of smooth
$\C$-valued k-forms and of (p,q)-forms on $X$. We will denote by
$\O_X(E), \A_X(E), \A_X^k(E), \A_X^{p,q}(E)$ the corresponding
sheaves of $E$-valued functions or forms.

Let $U\subset X$ be an open coordinate set with coordinates
$(z^1,\ldots,z^n)$. We assume that $E$ is trivial over $U$ and
we choose a holomorphic frame $\{ e_0,\ldots,e_r \}$ for $E$
over $U$. Relative to this frame vectors $v$ of $E$ can be
written
\begin{align*}
v= \sum_{i=0}^r v^i e_i
\end{align*}
and $E$-valued (p,q)-forms can be written
\begin{align*}
u= \sum_{i=0}^r \sum_{|I|=p,|J|=q} u^i_{IJ} dz^I \wedge 
                                      d\zb^J \tensor e_i
\end{align*}
where $dz^I= dz^{i_1}\wedge \ldots \wedge dz^{i_p}$ and
$J,I=\{ i_1,\ldots,i_p \}$ are increasing multiindices.

Let $h$ be a Hermitian metric on $E$. We represent it by a
Hermitian matrix
\begin{align*}
h=(h_{i\jb})_{0\le i,j\le r},\quad h_{i\jb}= <e_i, e_j>_h.
\end{align*}   
Associated to it there is the Chern connection $\nabla :\A(E)
\lra \A^1 (E)$ which can be represented by a matrix of 1-forms
\begin{align*}
(\th_i^j)_{0\le i,j\le r},\quad \nabla e_i = \sum_{j=0}^r 
\th_i^j \tensor e_j.
\end{align*}
The Chistoffel symbols of the first kind $\G^j_{ik}$ are
given by the relations
\begin{align*}
\th_i^j = \sum_{k=1}^n \G^j_{ik} dz^k.
\end{align*}
We have
\begin{align*}
\tag{3.1}
\th = (\de h)\cdot h^{-1},\quad \G^j_{ik}=\sum_{s=0}^r
\frac{\de h_{i\sb}}{\de z^k} \cdot h^{\sb j}.
\end{align*}
Here $(h^{\sb j})_{s,j}$ is the inverse of the matrix $h$.
The connection $\nabla$ induces unique $\C$-linear maps
$\nabla :\A^k(E) \lra \A^{k+1}(E)$ by enforcing the Leibnitz
rule: $\nabla (u \tensor v)= du \tensor v + (-1)^{|u|}u
\wedge \nabla v$. The composition $\nabla^2 :\A(E) \lra \A^2(E)$
is called the \emph{curvature} of $\nabla$ and has the remarkable
property that it is a tensor, i.e. it is a morphism of
$\A$-modules. $\nabla^2$ can be represented by a matrix
of (1,1)-forms
\begin{align*}
(\Th_i^j)_{0\le i,j\le r},\quad \nabla^2 e_i = \sum_{j=0}^r
\Th_i^j \tensor e_j.
\end{align*}
The Chistoffel symbols of the second kind $K^j_{ik\lb}$ are
given by the relations
\begin{align*}
\Th_i^j = \sum_{k,l=1}^n K^j_{ik\lb} dz^k \wedge d\zb^l.
\end{align*}
We have
\begin{align*}
\Th = d\th - \th \wedge \th = \deb \th,
\end{align*}
\begin{align*}
\tag{3.2}
K^j_{ik\bar{l}}= -\sum_{s=0}^r \frac{\de^2 h_{i\sb}}{\de z^k \de \zb^l}
h^{\sb j}+ \sum_{p,q,s=0}^r \frac{\de h_{i\bar{q}}}{\de z^k}\cdot
\frac{\de h_{p\bar{s}}}{\de \zb^l} h^{\bar{q}p} h^{\bar{s}j}.
\end{align*}
Given $v\in E_x$ we construct the (1,1)-form
\begin{align*}
\Th (v)=\frac{\i}{2\pi}
\frac{<\nabla^2 v, v >_h}{||v||^2_h}= \frac{\i}{2\pi} \frac{1}{||v||^2_h}
K^j_{i k \bar{l}} v^i h_{j \bar{s}} \bar{v}^s dz^k \wedge d \bar{z}^l.
\end{align*}
In the case of a line bundle $L$ the form $\Th (v)$ is nothing but the
first Chern form of the metric defined by
\begin{align*}
c_1 (L,h)= - dd^c \ \text{log}\ h.
\end{align*}
Given $\x \in$ T$_x X$ and $v\in E_x$ we define the \emph{holomorphic
bisectional curvature of $h$ along $\x$ and $v$} by
\begin{align*}
\tag{3.3}
k_x (\x,v)= \frac{<\nabla^2_{\x,\xb} v,v>_h}{||\x ||^2_{\f} 
||v||^2_h} = \frac{1}{||\x ||^2_{\f} ||v||^2_h}
\sum_{i,j=0}^r \sum_{k,l=1}^n K^j_{ik\lb}\x^k\xb^l v^i h_{j\sb}\vb^s.
\end{align*}

\noindent
{\bf (3.4) Definition:} Let $X$ be a complex manifold of dimension $n$
and $E$ a holomorphic vector bundle of rank $r+1\ge 2$ on $X$.
A \emph{Finsler metric} $h$ on $E$ is a function
\begin{align*}
h:E\lra [0,\infty)
\end{align*}
satisfying the following conditions:
\begin{enumerate}
\item[(i)] $h$ is continuous on $E$ and smooth on the complement of the
zero-section;
\item[(ii)] $h(\l v)= |\l |h(v)$ for $v\in E,\ \l \in \C$;
\item[(iii)] $h(v)>0$ if $v$ is non-zero;
\item[(iv)] $h_{| E_x \setminus \{ 0\}}$ is a strictly plurisubharmonic
function for all $x\in X$.\\
\end{enumerate}

The norm associated to a Hermitian metric is thus 
a particular case of Finsler metric.
We will see at (3.11) that a Finsler metric $h$ comes from a Hermitian
metric if and only if $h^2$ is of class $\mathcal{C}^2$ on $E$.

Let $\P(E)$ be the fiber bundle with fibers $\P(E_x)=$projective space of
lines through the origin in $E_x,\ \ x\in X$. We denote by
$\pi :\P(E)\lra X$ the projection onto the base. The tautological line
bundle $\L^{-1}= \O_{\P(E)}(-1)$ is the subbundle of $\p^*E$ whose fiber
at $(x,[v])$ consists of the line generated by $v$ inside $E_x$. The
bundle space of $\L^{-1}$ is the blow-up of $E$ along the zero-section;
let $\b :|\L^{-1}|\lra E$ be the blowing-up map. Outside the zero-sections
$\b$ is an isomorphism. To give $h$ satisfying properties (i), (ii) and
(iii) from above is equivalent to giving a Hermitian metric $\htil^{-1}$
on $\L^{-1}$ via the correspondence $\htil^{-1}= h\circ \b$.

In the sequel we will denote $E_o =E\setminus \{$zero-section\}.
We denote $p:E\lra X$ the projection onto the base and
by $q:E_o \lra \P(E)$ the quotient map. Thus
\begin{displaymath}
\xymatrix {
E_o \ar[rr]^{q} 
    \ar[dr]^{p} & & \P (E) \ar[dl]^{\p} \\
                & X
}
\end{displaymath}
is a commutative diagram. We will consider the function $G=h^2$
which is continuous on $E$ and smooth on $E_o$.

As before, let $U$ be a coordinate set which trivializes $E$.
On $p^{-1}U$ we have coordinates $(z^1,\ldots,z^n,v^0,\ldots,v^r)$.
We will consider the following smooth functions on $p^{-1}U\cap E_o$:
\begin{align*}
\tag{3.5}
G_{i\jb}=\frac{\de^2 G}{\de v^i \de \vb^j},\quad 0\le i,j\le r.
\end{align*}
By property (ii) for all $\l\in \C$ we have
\begin{align*}
\tag{3.6}
G(z,\l v)= |\l |^2 G(z,v).
\end{align*}
Differentiating with respect to $v^i$ and $\vb^j$ we obtain
\begin{align*}
\tag{3.7}
G_{i\jb}(z,\l v)= G_{i\jb} (z,v).
\end{align*}
Differentiating this relation with respect to $\l$ we obtain
\begin{align*}
\tag{3.8}
\sum_{p=0}^r v^p \frac{\de G_{i\jb}}{\de v^p} = 0 =
\sum_{q=0}^r \vb^q \frac{\de G_{i\jb}}{\de \vb^q}.
\end{align*}
Differentiating (3.6) with respect to $\l$ and $\bar{\l}$ we
obtain
\begin{align*}
\tag{3.9}
G(z,v)= \sum_{i,j=0}^r G_{i\jb} (z,\l v) v^i \vb^j
\end{align*}
which, in view of (3.7), takes the form
\begin{align*}
\tag{3.10}
G(z,v)= \sum_{i,j=0}^r G_{i\jb} (z,v) v^i \vb^j.
\end{align*}

\noindent 
If $h$ were a Hermitian metric with matrix $(h_{i\jb})$ we would
have $G= \sum_{i,j=0}^r h_{i\jb} v^i \vb^j$ and $G_{i\jb}(z,v)=
h_{i\jb}(z)$. Thus $h$ comes from a Hermitian metric if and only
if the functions $G_{i\jb}$ are constant along the fibers. In
fact, more can be said:

\noindent \\
{\bf (3.11) Remark:} $h$ comes from a Hermitian metric on $E$ if
and only if the function $G=h^2$ is of class $\mathcal{C}^2$ on $E$. Indeed,
if this were the case, we could take limit as $\l$ tends to 0
in (3.9) and we would obtain $G(z,v)= \sum_{i,j=0}^r G_{i\jb}(z,0)
v^i \vb^j$. Differentiating again with respect to $v^i$ and
$\vb^j$ we would get $G_{i\jb}(z,v)=G_{i\jb}(z,0)$, q.e.d.\\

Let us denote by $p_*:$T$E\lra $T$X$ the differential of $p$.
The kernel of $p_*$ is a holomorphic subbundle $\V\subset $T$E$
called the \emph{vertical tangent bundle}. Notice that, relative
to coordinates $(z^1,\ldots,z^n,v^0,\ldots,v^r)$ on $p^{-1}U$,
a frame of $\V$ is given by $\{ \frac{\de}{\de v^0},\ldots,
\frac{\de}{\de v^r}\}$. The transformation rule for this frame is
the same as for the frame $\{ e_0,\ldots, e_r\}$. This shows
that there is a canonical isomorphism $\V \isom p^*E$ identifying
$\frac{\de}{\de v^i}$ with $p^*e_i$. A vector $V$ of $\V$ can
be written
\begin{align*}
V=\sum_{i=0}^r V^i \frac{\de}{\de v^i}.
\end{align*}
A section of $\V$ of particular interest is the \emph{position
vector field}
\begin{align*}
P(z,v)=\sum_{i=0}^r v^i \frac{\de}{\de v^i}.
\end{align*}
Notice that $q^* \L^{-1}$ is the subbundle of $p^*E$ generated
by $P$.
Property (iv) tells us that the matrix $(G_{i\jb})_{0\le i,j\le r}$
is positive definite at each point of $E_o$. This allows us to
define a Hermitian metric on $\V_{|E_o}$ by putting
\begin{align*}
<V,W>_\V =\sum_{i,j=0}^r G_{i\jb}(z,v)V^i \overline{W}^j.
\end{align*}
This metric is nothing but pull-back of $h$ in the case when
$h$ is Hermitian. Notice that, by virtue of (3.10), we have
$G=||P||^2_\V$. Also notice that, by virtue of (3.9), $<,>_V$
is constant on the fibers of $q$ hence it descends to a
Hermitian metric of $\p ^*E$.

The Chern connection $\nabla^\V$ and its curvature make sense
just as for any Hermitian metric. $\nabla^\V$ has connection
matrix
\begin{align*}
\th_i^j =\sum_{k=1}^n \G^j_{ik} dz^k + \sum_{p=0}^r \g^j_{ip} dv^p
\end{align*}
where
\begin{align*}
\tag{3.12}
\G_{ik}^j=\sum_{s=0}^r \frac{\de G_{i\sb}}{\de z^k}\cdot G^{\sb j},
\quad \g^j_{ip}=\sum_{s=0}^r \frac{\de G_{i\sb}}{\de v^p}\cdot G^{\sb j}.
\end{align*}
The curvature matrix $\Th$ has now horizontal, vertical and
mixed components:
\begin{eqnarray*}
\Th_i^j & = & \sum_{k,l=1}^n K^j_{ik\bar{l}}\ dz^k\wedge d\zb^l +
             \sum_{p,q=0}^r \k^j_{ip\bar{q}}\ dv^p \wedge d\vb^q \\
        &   & + \sum_{k=1}^n \sum_{q=0}^r \m^j_{ik\bar{q}}\ dz^k
               \wedge d\vb^q + \sum_{l=1}^n 
                  \sum_{p=0}^r \n^j_{ip\bar{l}}\ dv^p \wedge d\zb^l,
\end{eqnarray*}
where
\begin{align*}
\tag{3.13}
K^j_{ik\bar{l}}= -\sum_{s=0}^r \frac{\de^2 G_{i\sb}}{\de z^k \de \zb^l}
G^{\sb j}+ \sum_{p,q,s=0}^r \frac{\de G_{i\bar{q}}}{\de z^k}\cdot
\frac{\de G_{p\bar{s}}}{\de \zb^l} G^{\bar{q}p} G^{\bar{s}j},
\end{align*}
\begin{align*}
\k^j_{ip\qb}= -\frac{\de \g^j_{ip}}{\de \vb^q},\ \m^j_{ik\qb}
=-\frac{\de \G^j_{ik}}{\de \vb^q},\ \n^j_{ip\lb}=
-\frac{\de \g^j_{ip}}{\de \zb^l}.
\end{align*}
Using (3.8) we obtain the following relations:
\begin{align*}
\tag{3.14}
\sum_{i=0}^r \g^j_{ip} v^i=0,
\end{align*}
\begin{align*}
\tag{3.15}
\sum_{i=0}^r \n^j_{ip\lb}v^i= \sum_{i=0}^r \k^j_{ip\qb}v^i=
\sum_{j,s=0}^r \m^j_{ik\qb} G_{j\sb} \vb^s =0.
\end{align*}

For any tangent vector $\z \in$ T$E_o$ of the form
\begin{align*}
\z =\sum_{k=1}^n a^k \frac{\de}{\de z^k} +
       \sum_{i=0}^r b^i \frac{\de}{\de v^i}
\end{align*}
we have
\begin{eqnarray*}
\nabla_\z^\V P & = & \sum_{i=0}^r (b^i +\sum_{j=0}^r \sum_{k=1}^n
\G^i_{jk} v^j a^k + \sum_{j,p=0}^r \g^i_{jp} v^j b^p)
\frac{\de}{\de v^i}\\
& = & \sum_{i=0}^r (b^i +\sum_{j=0}^r \sum_{k=1}^n
\G^i_{jk} v^j a^k )
\frac{\de}{\de v^i} \quad \quad \text{by (3.14)}.
\end{eqnarray*}
This calculation shows that the linear map of bundles
\begin{align*}
\nabla^\V P:\text{T}E_o \lra \V,\quad \z \lra \nabla_\z^\V P
\end{align*}
is a surjection. Thus the kernel of $\nabla^\V P$ is a smooth
subbundle $\H \subset $ T$E_o$ which we call the \emph{horizontal
tangent bundle}. We have the smooth decomposition
\begin{align*}
\text{T}E_o =\H \oplus \V
\end{align*}
and at any point $(x,v)\in E_o$ the differential $p_*:\H_{(x,v)}
\lra $ T$_x X$ is an isomorphism. Under this isomorphism a vector
$\x \in $ T$_x X$ corresponds to a vector $\x^\H \in \H_{(x,v)}$
called the \emph{horizontal lift of} $\x$. The vertical tangent
fields together with the horizontal lifts of the tangent fields
$\de_k =\de /\de z^k$ give a smooth frame of T$E_o$:
\begin{align*}
\tag{3.16}
\left\{
\begin{array}{ll}
\de / \de v^i, & 0\le i\le r, \\
\de_k^\H = \de / \de z^k - \sum_{i,j=0}^r \G^i_{jk} v^j
\de / \de v^i, & 1\le k\le n.
\end{array}
\right.
\end{align*}
The dual basis is
\begin{align*}
\tag{3.17}
\left\{
\begin{array}{ll}
dz^k, & 1\le k\le n,\\
\z^i = dv^i + \sum_{k=1}^n \sum_{j=0}^r \G^i_{jk} v^j dz^k,
& 0\le i\le r.
\end{array}
\right.
\end{align*}

Just as in the Hermitian case, for any point $(x,v)\in E_o$
there is a holomorphic frame of $E$ on a neighbourhood of $x$
with respect to which we have
\begin{align*}
\tag{3.18}
G_{i\jb} (x,v)= \d_{ij},\quad \frac{\de G_{ij}}{\de z^k}(x,v)=
0,\quad 0\le i,j\le r,\ 1\le k\le n.
\end{align*}
Here $\d_{ij}$ is the Kronecker symbol. We call such a frame
\emph{normal at} $(x,v)$. Relative to such a frame we have
$\de_k^\H =\de /\de z^k,\ \z^i=dv^i$ at $(x,v)$.

We claim that for any non-zero vector $v \in E_x$ the following
expression gives us a well-defined (1,1)-form on T$_x X$:
\begin{align*}
\tag{3.19}
<K(\cdot,\cdot)v,v>_\V = \sum_{i,j,p=0}^r \sum_{k,l=1}^n
K^p_{ik\lb} v^i G_{p\jb} \vb^j dz^k\wedge d \zb^l.
\end{align*}
To see this we will show that for any $\x_1,\x_2 \in$ T$_x X$
we have
\begin{align*}
\tag{3.20}
<K(\x_1,\bar{\x}_2)v,v>_\V= <\nabla^2_{\x_1^\H,\bar{\x}_2^\H}
P,P>_\V
\end{align*}
where $\nabla^2$ is the curvature of $\nabla^\V$. Indeed,
adopting Einstein's summation convention, we have
\begin{eqnarray*}
<\nabla^2 P,P>_\V & = & K^j_{ik\lb} v^i G_{j\sb} \vb^s
d z^k \wedge d\zb^l + \k^j_{ip\qb} v^i G_{j\sb} \vb^s dv^p
\wedge d\vb^q \\
& & +\m^j_{ik\qb} v^i G_{j\sb}\vb^s dz^k\wedge d\vb^q +
       \n^j_{ip\lb} v^i G_{j\sb}\vb^s dv^p\wedge d\zb^l.
\end{eqnarray*}
By (3.15) the last three terms on the right-hand-side vanish,
yielding (3.20). As a byproduct we get the relation
\begin{align*}
\Th(P)= \frac{\i}{2\pi}\cdot \frac{1}{||P||^2_\V}
\sum_{i,j,s=0}^r K^j_{ik\lb} v^i G_{j\sb} \vb^s dz^k \wedge
d\zb^l.
\end{align*}

In order to define the holomorphic bisectional curvature of $h$
we need to first put a metric on T$E_o$. This is done by making
the decomposition T$E_o=\H\oplus \V$ an orthogonal decomposition.
Namely, any tangent vector $\z \in $ T$E_o$ can be written
$\z= \x^\H + V$ and we put $||\z||^2=||\x||^2_\f + ||v||^2_\V$.
We then define the \emph{holomorphic bisectional curvature
of $h$ along $\z$ and $v$} just as at (3.3):
\begin{align*}
\tag{3.21}
k(\z,v)= \frac{<\nabla^2_{\z,\bar{\z}}v,v>_\V}{||\z ||^2 ||v||^2_h} 
\qquad \qquad \qquad \qquad \qquad \qquad \quad
\end{align*}
\begin{eqnarray*}
\qquad  \qquad \quad \quad \
& = & \frac{1}{||\z ||^2 ||v||^2_h} \sum_{i,j,s=0}^r
                 \Big\{ \sum_{k,l=1}^n K^j_{ik\lb} a^k \ab^l +
                    \sum_{p,q=0}^r \k^j_{ip\qb} b^p \bb^q \\
&   & \qquad \quad + \sum_{k=1}^n \sum_{q=0}^r \m^j_{ik\qb} a^k \bb^q +
                 \sum_{l=1}^n \sum_{p=0}^r \n^j_{ip\lb} b^p \ab^l
                       \Big\} v^i G_{j\sb}\vb^s.
\end{eqnarray*}
The knowledge that (3.19) is well defined allows us to
introduce the notion of \emph{horizontal holomorphic
bisectional curvature of $h$ along $\x\in $ T$_x X$ and $v\in E_x$}:
\begin{align*}
\tag{3.22}
k_x (\x,v)= k_{(x,v)}(\x^\H,v)=
\frac{<\nabla^2_{\x^\H,\bar{\x}^\H}P,P>_\V}
{||\x ||^2_\f ||v||^2_h} 
\end{align*}
\begin{align*}
\qquad \qquad \quad \ \ = \frac{1}{||\x ||^2_\f ||v||^2_h}
\sum_{i,j,s=0}^r \sum_{k,l=1}^n 
K^j_{ik\lb} \x^k \bar{\x}^l v^i G_{j\sb} \vb^s.
\end{align*}

In the remaining part of this section we will compute
$c_1 (\L,\htil)$. Recall that, under the isomorphism
$\V \isom p^*E$, we can identify $q^* \L^{-1}$
with the sub-bundle of $\V$ spanned by $P$.
Relation (3.10) tells us that the pull-back metric
$q^* \htil^{-1}$ is nothing but the induced metric from $\V$.
Therefore we have
\begin{align*}
q^* c_1(\L,\htil)= dd^c \text{log} ||P||^2_\V.
\end{align*}
We claim that
\begin{align*}
dd^c \text{log} ||P||^2_\V & = &
\frac{\i}{2\pi}\cdot \frac{1}{G^2} \Big\{ \sum_{i,j=0}^r G G_{i\jb}
\z^i \wedge \bar{\z}^j - 
\sum_{i,j,p,q=0}^r G_{i\qb} G_{p\jb} \vb^q v^p
\z^i \wedge \bar{\z}^j \Big\}
\end{align*}
\begin{align*}
\tag{3.23}
\quad \quad -\frac{\i}{2\pi} \cdot \frac{1}{G}
\sum_{i,j,s=0}^r \sum_{k,l=1}^r
K^j_{ik\lb} v^i G_{j\sb}\vb^s dz^k \wedge d\zb^l.
\end{align*}
The part involving terms $\z^i \wedge \bar{\z}^j$ will be called
the \emph{vertical component of} $c_1 (\L,\htil)$ and will be
denoted $c_1 (\L,\htil)^\V$. The part involving terms
$dz^k \wedge d\zb^l$ will be called the \emph{horizontal
component of} $c_1 (\L,\htil)$ and will be denoted
$c_1 (\L,\htil)^\H$. Notice that $c_1 (\L,\htil)^\V$ is
semi-positive definite and its restriction to the fibers
of $\P(E)$ is positive-definite. Indeed, relative to a normal
frame for $E$ at $(x,v)$ we have
\begin{align*}
c_1 (\L,\htil)^\V = \frac{\i}{2\pi}\cdot \frac{1}{||v||^2}
\sum_{i=0}^r dv^i \wedge d\vb^i
- \frac{\i}{2\pi}\cdot \frac{1}{||v||^4}
\sum_{i,j=0}^r \vb^i v^j dv^i \wedge d\vb^j
\end{align*}
with $||v||^2= \sum_{i=0}^r |v_i |^2$. The above is nothing
but the Fubini-Study form on projective space which is known
to be positive-definite. Notice that
\begin{align*}
c_1 (\L,\htil)^\H = - \Th (P).
\end{align*}
In normal coordinates at $(x,v)$ we have
\begin{align*}
\tag{3.24}
c_1 (\L,\htil)= \o - \Th(P)
\end{align*}
where $\o$ is the Fubini-Study form on the first component of
\begin{align*}
\text{T}_{(x,[v])}\P(E)= \text{T}_{[v]} \P(E_x)
\oplus \text{T}_x X.
\end{align*}
Thus $\L$ equipped with $\htil$ is a positive line bundle
if and only if the horizontal holomorphic bisectional
curvature of the Finsler metric $h$ is negative! This has
an important consequence: assume $X$ to be compact; then $E^*$
is ample if and only if there is a Finsler metric (with
only properties (i), (ii) and (iii)) on $E$ having negative
horizontal holomorphic bisectional curvature.
We will not need this fact, but it
is worth mentioning because it illustrates the philosophy
we advertized in the introduction: geometric properties
of $E^*$ should be defined by means of $\O_{\P(E)}(1)$. \\ \\

%
%

\section{Hermitian metrics of finite order}

Let $E$ be a holomorphic vector bundle on the special affine 
variety $X$. We recall from \cite{griffiths-cornalba}
the notion of Hermitian metric on $E$ 
of finite order. It is defined
by means of an estimate on the holomorphic bisectional
curvature. The theorem of Griffiths and Cornalba (4.4) states
that a holomorphic line bundle on $X$
admits a unique finite order structure. In fact, as we
explained in the introduction, Griffiths and Cornalba
gave four definitions for the notion of
finite order vector bundle which are equivalent in the
case of line bundles.
One of the crucial steps in the proof of the four
equivalences is their vanishing theorem (4.6) for the
sheaf of sections of finite order.

\noindent \\
{\bf (4.1) Definition:} A Hermitian holomorphic vector
bundle $(E,h)$ on $X$ is said to \emph{have order $\l$}
if its holomorphic bisectional curvature is of order $\r^\l$:
there is $\k >0$ such that
\begin{align*}
|k_x (\x,v)|\le \k \r^\l
\end{align*}
for all $x\in X,\ \x \in $ T$_x X,\ v\in E_x$.

Notice that, adopting the notation following (3.2),
the above estimate is equivalent to the
inequality
\begin{align*}
|\Th (v)|\le \k \r^\l \f.
\end{align*}
Also notice that a line bundle $(L,h)$ has finite order if and only if
\begin{align*}
|c_1 (L,h)| \le \k \r^\l \f.
\end{align*}

\noindent
{\bf (4.2) Remark:} The above definition is preserved under
standard operations with vector bundles: direct sum, tensor
product, dualization, symmetric power, exterior power etc.
In particular, if $E_1$ and $E_2$ are Hermitian vector bundles
of finite order then so is ${\mathcal Hom}(E_1, E_2)
\isom E_1^* \tensor E_2$. \\

Sections of finite order of a Hermitian vector bundle
are defined the same way as functions: for a global section
$\s$ of $E$ we introduce the \emph{maximum modulus function}
\begin{align*}
M_\s (r)= \text{log max}\{ |\s(x)|_h,\ \r (x)\le r\}.
\end{align*}
We say that $\s$ \emph{has order $\l$ in the sup-norm sense}
if $M_\s (r)=O(r^\l)$. In view of the previous remark this
gives us a notion of \emph{finite order morphism} between
two Hermitian vector bundles. The composite of morphisms of
finite order is of finite order, too. Thus we have constructed
the category of \emph{Hermitian holomorphic vector bundles
of finite order}. In this category a vector bundle of
rank $r$ is trivial, i.e. it is
isomorphic to the trivial vector bundle
$\underline{\C}^r$, if and only if it posesses $r$ global
sections of finite order which generate the fiber
at every point. We will denote by Vect$^r_{\fo}(X)$ the set
of all Hermitian holomorphic vector bundles of rank $r$
on $X$ modulo isomorphisms in the finite order category.
Vect$^1_{\fo}(X)$ together with the
tensor product forms a group which we call the
\emph{Picard group of finite order line bundles on $X$}
and denote Pic$_{\fo}(X)$.
As we explained in the introduction, the motivation for
our work is the following:

\noindent \\
{\bf (4.3) Question:} Is the canonical map Vect$^r_{\fo}(X)
\lra $ Vect$^r_{\text{hol}} (X)$ a bijection? \\

The question was answered in the affirmative for the case of
line bundles.

\noindent \\
{\bf (4.4) Theorem} (Griffiths and Cornalba 
\cite{griffiths-cornalba}):
\emph{The canonical map
\begin{align*}
\text{Pic}_{\fo} (X)\lra \text{Pic}_{\text{hol}}(X)
\end{align*}
is an isomorphism of groups. In other words, every holomorphic
line bundle on $X$ admits a metric of finite order and a
Hermitian holomorphic line bundle of finite order is
trivial in the holomorphic category if and only if it is
trivial in the finite order category.}

\noindent \\
{\bf (4.5) Definition:} Let $(E,h)$ be a holomorphic Hermitian
vector bundle. A global section $\s$ of $E$ is said to
\emph{have order $\l$} if there is $\k \ge 0$ such that
\begin{align*}
\int_X |\s |^2_h \cdot e^{-\k \r^\l} \Phi \ < \ \infty.
\end{align*}
Here $\Phi =\f^n$ is the volume form of $\f$. We define the
\emph{sheaf $\O_\l (E)$ of germs of holomorphic sections
of $E$ of order $\l$} as follows: $\O_\l (E)$ is a sheaf on
$\Xb$; for each open set $U\subset \Xb$ the space of sections
$\O_\l (U,E)$ consists of those holomorphic sections
$\s$ of $E$ defined on $U \cap X$ with the property
that around each point at infinity $x\in (\Xb \setminus X)\cap
U$ there is a neighbourhood $W\subset U$ and a constant
$\k >0$, depending on $x$ and $\s$, such that 
\begin{align*}
\int_{W\cap X} |\s |^2_h \cdot e^{-\k \r^\l} \Phi \ < \ \infty.
\end{align*}
Similarly, we define the \emph{sheaf $\O_{\fo}(E)$ of
germs of holomorphic sections of $E$ of finite order}
by the same estimate as above with the additional requirement
that $\l$ depend on $\s$ and on the point at infinity $x$.

It is easily seen that $\O_\l (E)$ and $\O_{\fo}(E)$ are
modules over $\O_\l$, respectively over $\O_{\fo}$. Their
restrictions to $X$ coincide with the sheaf $\O (E)$ of
germs of sections of $E$ because in (4.5) there are no
conditions at the points $x\in X$ away from infinity.
By the compactness of the divisor at infinity it is also
transparent that the spaces of global sections
$\O_\l (\Xb,E)$ and $\O_{\fo} (\Xb,E)$ are nothing but
the spaces of global sections of $E$ of order $\l$,
respectively of finite order.\\

Let us now assume that $E$ has order $\l$. Then a global
section $\s$ of $E$ has finite order in the sup-norm sense
if and only if it has finite order in the sense of (4.5).
One direction is immediate: if $M_\s (r)= O(r^\l)$ then
$|\s |_h \le e^{\k \r^\l}$ for some $\k \ge 0$ and the
integral
\begin{align*}
\int_X |\s |^2_h \cdot e^{-2 \k \r^\l} \Phi
\end{align*}
is finite because $\Phi$ has finite volume. For the other
direction we refer to \cite{griffiths-cornalba}.

\noindent \\
{\bf (4.6) Theorem:}
\emph{Let $E$ be a Hermitian holomorphic
vector bundle of order $\l$ on a special affine variety
$X$. Then, for all $q\ge 1$, we have}
\begin{align*}
H^q (\Xb, \O_\l(E))=0,\quad H^q (\Xb, \O_{\fo}(E))=0.
\end{align*}

We now prepare to move from the Hermitian to the Finsler
case. We begin by translating definition (4.5) ``upstairs''
on $\P (E)$. There is a canonical isomorphism of vector
bundles $\pi_* (\L) \isom E^*$ which associates to a
section $\s$ of $E^*$ the section $\st$ of $\L$ defined
by
\begin{align*}
\tag{4.7}
\st_{(x,[v])} (t v)= \s (tv),
\end{align*}
for all $x\in X,\ v\in E_x \setminus \{ 0\},\ t\in \C$.
Recall that the Hermitian metric $h$ on $E$ (which in
particular is a Finsler metric) induces a metric $\htil$
on $\L$. Denote by $h^*$ the dual metric on $E^*$.
It happens that
\begin{align*}
\tag{4.8}
|\s (x)|_{h^*}= \text{max} \{ |\st (x,[v])|_{\htil},\ \ 
v\in E_x \setminus \{ 0\} \}.
\end{align*}

Assume that $E$ has a Hermitian metric $h$ of order $\l$.
We claim that there is a positive constant $\k$ such that the expression
\begin{align*}
\tag{4.9}
\ft = c_1 (\L,\htil) + dd^c (\k \r^\l)
\end{align*}
gives a positive (1,1)-form on $\P(E)$. We check this at
a point $(x,[v])$. Since positivity is preserved under
change of coordinates we can work with a frame of $E$
which is normal at $x$. By our hypothesis (4.1) we have
$|\Th (v)|\le \k \r^\l \f$. We also have
\begin{eqnarray*}
\frac{\k \l}{2} \r^\l \f  =
\frac{\k \l}{2} \r^\l dd^c \t   & \le & \frac{\k \l}{2} \r^\l dd^c \t
               + \frac{\k \l^2}{4} \r^\l d\t \wedge d^c\t \\
                         & =   & \k \ dd^c e^{\l \t /2}     \\
                         & =   & dd^c (\k \r^\l)  \\
                         & \le & \frac{\k \l}{2} \r^\l dd^c \t
              + \frac{\k \l^2}{4} \r^{\l+c}\f \quad \quad \qquad
                                            \text{by (1.8)}.     
\end{eqnarray*}
From this calculation we get the estimate
\begin{align*}
\tag{4.10}
\frac{\k \l}{2}\r^\l \f \le -\Th(v) + dd^c (\k \r^\l)
\le 2\k \r^{\l +c}\f.
\end{align*}
Combining (3.24) with (4.10)
we get the following estimate at $(x,[v])$:
\begin{align*}
\tag{4.11}
\o + \frac{\k \l}{2} \r^\l \f \le \ft \le \o + 2 \k \r^{\l +c}\f.
\end{align*}
This insures the positivity of $\ft$. In the sequel $\P(E)$
will be considered equipped with the K\"ahler metric induced
by $\ft$. Since $\r^\l \circ \pi$ is a plurisubharmonic
exhaustive function, we can arrange, possibly by choosing
a larger $\k$, that this metric on $\P(E)$ be complete.

\noindent \\
{\bf (4.12) Claim:}
\emph{Assume that $E$ is a holomorphic Hermitian
vector bundle of order $\l$. Let $U\subset X$ be an open subset
and $\s$ a section of $E^*$ over $U$. Then there exists
$\k_1 \ge 0$ such that}
\begin{align*}
\int_U |\s |^2_{h^*} \cdot e^{-\k_1 \r^\l} \Phi \ < \ \infty
\end{align*}
\emph{if and only if there exists $\k_2 \ge 0$ such that}
\begin{align*}
\int_{\pi^{-1}(U)} |\st |^2_{\htil}\cdot e^{-\k_2 \r^\l}\Phit
                                                \ < \ \infty.
\end{align*}

\noindent
\emph{Proof:} Let us denote by $\ft^\V$ and $\ft^\H$ the
vertical, respectively the horizontal part of $\ft$, as defined
at the end of \S 3. From (4.10) we have
\begin{align*}
\tag{4.13}
\frac{\k \l}{2} \r^\l \f \le \ft^\H \le 2\k \r^{\l +c}\f
\end{align*}
where, by an abuse of notation, we write $\f$ instead of
$\pi^* \f$. We have
\begin{align*}
\Phit = \ft^{n+r} = {n+r \choose n}
(\ft^\H)^n \wedge (\ft^\V)^r
\end{align*}
which, combined with (4.13), gives
\begin{align*}
\tag{4.14}
C^{-1} \r^{n\l} \f^n \wedge (\ft^\V)^r \le 
\Phit \le 
C \r^{n\l +nc} \f^n \wedge (\ft^\V)^r
\end{align*}
for a fixed positive constant $C$. 
This estimate is useful because one can apply Fubini's theorem
to $|\st |^2_{\htil} \cdot e^{-\k \r^\l} \f^n \wedge (\ft^\V)^r$.
More precisely, one can integrate this form first vertically along
the fibers of $\P(E)$ and then horizontally along $X$.
This finishes the proof of the claim because there
is a positive constant $A$ depending
only on $r$ such that for all $x\in X$
\begin{align*}
\tag{4.15}
\int_{[v]\in \P(E_x)} |\st (x,[v])|^2_{\htil} (\ft_{|\P(E_x)})^r
= A |\s (x)|^2_{h^*}.
\end{align*}
\\

%
%

\section{Finsler Metrics of Finite Order. The Vanishing Theorem}

We explained in the introduction
the difficulty one encounters in trying to generalize
(4.4) and the four equivalences to bundles 
of rank greater than 1. It seems to us
that the correct antidote is to translate the definition
of finite order ``upstairs'' on $\P(E)$. So we define
Finsler metrics of finite order by a very similar
estimate on the holomorphic bisectional curvature.
We then define sections of finite order and we prove 
that they span at every point, cf. (5.23). This
we achieve by means of the vanishing theorem (5.12)
which generalizes (4.6). 

Traditionally, vanishing
theorems are given either on Stein manifolds
or on compact K\"ahler manifolds. Our case here is a
hybrid: we will have to solve the $\deb$-equation on $\P(E)$
whose base is Stein (because it is affine)
while its fibers are compact.

The original proof of Kodaira's Vanishing Theorem makes
use of Hodge theory: one knows that on a compact
K\"ahler manifold cohomology classes can be represented
by harmonic forms and one argues, using the a priori
estimate, that such forms do not exist.
For non-compact manifolds this argument does not
work because we do not know if the Hodge representation
theorem holds. Instead we will use a very potent technique
developed by H\"ormander in \cite{hormander-estimates}
involving some rudiments of functional analysis.
We synthetize this technique in proposition (5.3)
and theorem (5.4) from below.
We refer to \cite{hormander}, \cite{hormander-estimates},
\cite{andreotti-vesentini} and \cite{kodaira}.

Let $Y$ be a complex manifold of dimension $m$ equipped with
a positive K\"ahler form $\o$ which induces a complete
metric. Relative to local holomorphic coordinates
$(z^1,\ldots,z^m)$ we write
\begin{align*}
\o =\i \sum_{i,j=1}^m g_{i \jb}\ dz^i \wedge d\zb^j.
\end{align*}
The condition that $\o$ be positive definite means that
the matrix $(g_{i\jb})_{1\le i,j\le m}$ is Hermitian and
positive-definite. The condition that $\o$ be K\"ahler
means that $d \o =0$. This is a very natural condition
to consider
because it is equivalent to saying that the complex
structure on the real tangent bundle of $Y$
is compatible with the Levi-Civita connection associated
to the induced Riemannian metric on the real tangent
bundle of $Y$.

Let Ric$(\o)$ be the Ricci curvature of $\o$. It is a
(1,1)-form given in local coordinates by
\begin{align*}
\text{Ric}(\o)= dd^c \text{log(det}(g_{i \jb})).
\end{align*}
It is nothing but the first Chern form of the canonical
line bundle $K_Y := \wedge^m $ T$^* Y$ equipped with
the metric induced by $\o$.

Recall from \S 3
that the first Chern form of a holomorphic line
bundle $L$ equipped with a Hermitian metric $h$ is given
by
\begin{align*}
c_1 (\L,h)= -dd^c \text{log}(h).
\end{align*}

Let us choose an orthonormal frame $\{ \x_1, \ldots , \x_m \}$
for the tangent space at a point $y\in Y$ and a unitary
vector $e\in L_y$. Relative to the dual frame
$\{ d\x^1, \ldots, d\x^m \}$ for T$^*_y Y$ we write
\begin{align*}
\text{Ric}(\o)= \i \sum_{i,j=1}^m
R_{i \jb}\ d\x^i \wedge d\xb^j,\\
c_1(\L)= \i \sum_{i,j=1}^m K_{i \jb}\ d\x^i \wedge d\xb^j.
\end{align*}
For an $L$-valued (0,q)-form
\begin{align*}
u= \sum_{|I|=q} u_I \ d\xb^I \tensor e
\end{align*}
we define the pointwise operators
\begin{align*}
<Ru,u>= q \sum_{i,j=1}^m \sum_{|I|= q-1} R_{i \jb}\ u_{iI} \bar{u}_{jI},
\end{align*}
\begin{align*}
<Ku,u>= q \sum_{i,j=1}^m \sum_{|I|= q-1} K_{i \jb}\ u_{iI} \bar{u}_{jI},
\end{align*}
and their integrated versions
\begin{align*}
(Ru,u)= \int_Y <Ru,u> dV,
\end{align*}
\begin{align*}
(Ku,u)= \int_Y <Ku,u> dV.
\end{align*}
Here $u$ is assumed to have compact support and $dV=\o^m$ is the
volume form of $\o$. In the sequel we will denote by
$\DD^{p,q}(Y,L)$ the space of smooth $L$-valued (p,q)-forms
with compact support. The Hermitian inner products on
T$^*Y$ and on $L$ induce a Hermitian inner product $<,>$
on $\wedge^{p,q}$T$^*Y \tensor L$. Given $u,v\in \DD^{p,q}(Y,L)$
we put
\begin{align*}
(u,v)= \int_Y <u,v>\ dV.
\end{align*}
Clearly $(\cdot,\cdot)$ defines a Hermitian inner product
on $\DD^{p,q}(Y,L)$. Its associate norm is given by
$||u||^2 = (u,u)$. The operator
\begin{align*}
\deb : \DD^{p,q}(Y,L) \lra \DD^{p,q+1}(Y,L)
\end{align*}
has a formal adjoint
\begin{align*}
\del : \DD^{p,q+1}(Y,L) \lra \DD^{p,q}(Y,L)
\end{align*}
given by the condition $(\deb u, v)=(u,\del v)$ for all
$u\in \DD^{p,q}(Y,L),\ v\in \DD^{p,q+1}(Y,L)$. We mention in
passing that $\del =-\star D' \star$ where $\star$ is the
Hodge-$\star$ operator while $D'$ is the (1,0)-component
of the Chern connection of $L$.

Our key ingredient towards proving vanishing theorems is
the following inequality which we state in a particular case
that is of interest to us. See \cite{kodaira} p. 124 for
the full statement and p. 126 for the statement from below.

\noindent \\
{\bf (5.1) Weitzenb\"ock Inequality:}
\emph{Let $Y$ be a K\"ahler
manifold. Let $L$ be a holomorphic Hermitian line bundle
on $Y$ and $u$ a smooth $L$-valued (0,q)-form with compact support.
Then we have:}
\begin{align*}
||\deb u||^2 +||\del u||^2 \ge (Ku,u) - (Ru,u).
\end{align*}

\noindent \\
{\bf (5.2) Definition:} We say that $L$ is \emph{(p,q)-elliptic}
if there is a positive constant $\e$, called \emph{ellipticity
constant}, such that for all $u\in \DD^{p,q}(Y,L)$ we have the
estimate
\begin{align*}
||\deb u||^2 +||\del u||^2 \ \ge \ \e ||u||^2.
\end{align*}

\noindent \\
{\bf (5.3) Proposition:}
\emph{Assume $L$ is (p,q)-elliptic
and that the metric on $Y$ is complete and K\"ahler. 
Then for any $u\in $ L$^{p,q}(Y,L)$ with $\deb u=0$ there
is $v\in $ L$^{p,q-1}(Y,L)$ with $\deb v=u$. In addition
$v$ is smooth.} \\

Here $\text{L}^{p,q}(Y,L)$ is the completion of
$\DD^{p,q}(Y,L)$ relative to the norm $||u||^2$.
Notice here the crucial assumption that $Y$ be complete.
An equivalent condition is that any ball in $Y$ be
relatively compact. We need this for a density argument
to make sure that (5.1) holds not only for smooth forms
but also for square-integrable forms.
See \cite{andreotti-vesentini}, p. 92, lemma 4.

\noindent \\
{\bf (5.4) Theorem:}
\emph{Let $Y$ be a complex manifold
equipped with a complete K\"ahler metric. Let $L$ be a
holomorphic line bundle on $Y$ equipped with a Hermitian
metric $h$. We denote by $\o$ the K\"ahler form of $Y$.
Assume that there is a positive constant $\e$ such that}
\begin{align*}
c_1 (L,h) - \text{Ric}(\o) \ \ge \ \e \ \o.
\end{align*}
\emph{Let $q$ be a positive integer.
Then for any $u\in $ L$^{0,q}(Y,L)$ with $\deb u=0$ there
is $v\in $ L$^{0,q-1}(Y,L)$ with $\deb v=u$. In addition
$v$ is smooth.}

\noindent \\
\emph{Proof:} By hypothesis for any $u\in \DD^{0,q}(Y,L)$
we have
\begin{align*}
(Ku,u)-(Ru,u) \ \ge \ q\e ||u||^2.
\end{align*}
Combining this with the Weitzenb\"ock Inequality we
obtain
\begin{align*}
||\deb u||^2 +||\del u||^2 \ \ge \ q\e ||u||^2.
\end{align*}
Hence $L$ is (0,q)-elliptic with ellipticity constant
$q\e$. The theorem now follows from (5.3).

\noindent \\
{\bf (5.5) Definition:} Let $E$ be a holomorphic vector
bundle on $X$ and $h$ a Finsler metric on $E$. We say that
\emph{$h$ has order $\l$} if its holomorphic bisectional
curvature is of order $\r^\l$: there is $\k>0$ such that
\begin{align*}
|k (\z,v)|\le \k \r^\l
\end{align*}
for all $y\in E_o,\ \z \in \text{T}_y E_o,\ v\in E_x$
where $p (y)=x$.

Notice that this definition encompasses (4.1). For a line
bundle a Finsler metric is the same thing as a Hermition metric,
so (5.5) is relevant only when the rank of $E$ is at least 2,
which we will assume henceforth.

Notice that if we take $\z$ to be the horizontal lift of a tangent
vector from $X$ and $v$ to be the position vector field we obtain
that the horizontal holomorphic bisectional curvature is also
of order $\r^\l$:
\begin{align*}
\tag{5.6}
| k_x (\x,v)|\le \k \r^\l
\end{align*}
for all $x\in X,\ \x \in \text{T}_x X,\ v\in E_x$. Adopting the
notation below (3.20) the above estimate is equivalent to the
inequality
\begin{align*}
|\Th (P)|\le \k \r^\l \f.
\end{align*}
In view of (3.24) and using the same argument as in the Hermitian
case we can find a positive constant $\k$ such that the expression
\begin{align*}
\tag{5.7}
\ft = c_1 (\L,\htil)+dd^c (\k \r^\l)
\end{align*}
defines a K\"ahler form on $\P(E)$ which induces a complete metric.
Estimate (4.13) also holds:
\begin{align*}
\tag{5.8}
\frac{\k \l}{2} \r^\l \f \le \ft^\H \le 2\k \r^{\l+c}\f.
\end{align*}
Moreover, the constant $\k$ can be so chosen that
\begin{align*}
|c_1 (G_{i\jb})| \le \k \r^\l \ft.
\end{align*}
This is so because an estimate on the holomorphic bisectional
curvature of a vector bundle $E$ gives an estimate
of the curvature of the induced metric on the determinant of $E$.
From the proof of (5.12) it transpires that the Ricci curvature of
$\P(E)$ is, in essence, equal to $c_1 (G_{i\jb})$ plus some other
terms that are controllable. Thus, roughly speaking, the above estimate
says that the Ricci curvature of $\P(E)$ has polynomial growth.
Taking horizontal and vertical parts we obtain the following inequalities:
\begin{align*}
\tag{5.9}
|c_1 (G_{i\jb})^\H|\le \k \r^\l \f,
\end{align*}
\begin{align*}
\tag{5.10}
|c_1(G_{i\jb})^\V|\le \k \r^\l c_1 (\L,\htil)^\V.
\end{align*}
In fact, (5.6), (5.9) and (5.10) are all that we will need in the sequel.
It is not clear to us if these estimates imply (5.5), in other words
it seems to us that (5.5) is stronger than (5.6), (5.9) and (5.10).
We could have chosen the three latter estimates as our definition
of finite order Finsler metric, but then the similarity with
the Hermitian case would have been obscured.

Our Finsler metric is on $E$ but we will be concerned with
and we will prove a vanishing theorem for the sections of
order $\l$ of $E^*$. We define the latter in the spirit of
(4.12), by using the one-to-one correspondence (4.7) between
sections of $E^*$ and sections of $\L$. Before doing that
we fix a Hermitian metric $g$ of finite order on det$(E)$.
Such a metric exists by (4.4). Under the isomorphism
\begin{align*}
\L \isom \L \tensor \p^* (\text{det}(E^*))\tensor
\p^* (\text{det}(E))
\end{align*}
we consider the following metric on $\L$:
\begin{align*}
\ltil = \htil \cdot \text{det}(G_{i\jb})^{-1} \cdot \p^*(g).
\end{align*}

\noindent
{\bf (5.11) Definition:} Let $E$ be a holomorphic vector
bundle on $X$ equipped with a Finsler
metric of order $\l$.
A global section $\s$ of $E^*$ is
said to \emph{have order} $\l$ if there is $\k >0$
such that
\begin{align*}
\int_{\P(E)} |\st |^2_{\ltil} \cdot e^{-\k \r^\l} \Phit
\ < \ \infty.
\end{align*}

We define the
\emph{sheaf $\O_\l (E^*)$ of germs of holomorphic sections
of $E^*$ of order $\l$} as follows: $\O_\l (E^*)$ is a sheaf on
$\Xb$; for each open set $U\subset \Xb$ the space of sections
$\O_\l (U,E^*)$ consists of those holomorphic sections
$\s$ of $E^*$ defined on $U \cap X$ with the property
that around each point at infinity $x\in (\Xb \setminus X)\cap
U$ there is a neighbourhood $W\subset U$ and a constant
$\k >0$, depending on $X$ and $\s$, such that
\begin{align*}
\int_{\p^{-1}(W\cap X)}
|\st |^2_{\ltil} \cdot e^{-\k \r^\l} \Phit \ < \ \infty.
\end{align*}

Similarly, we define the \emph{sheaf $\O_{\fo}(E^*)$ of
germs of holomorphic sections of $E^*$ of finite order}
by the same estimate as above with the additional requirement
that $\l$ depend on $\s$ and on the point at infinity $x$.

It is easily seen that $\O_\l (E^*)$ and $\O_{\fo}(E^*)$ are
modules over $\O_\l$, respectively over $\O_{\fo}$. Their
restrictions to $X$ coincide with the sheaf $\O (E^*)$ of
germs of sections of $E^*$ because in (5.11) there are no
conditions at the points $x\in X$ away from infinity.
By the compactness of the divisor at infinity it is also
transparent that the spaces of global sections
$\O_\l (\Xb,E^*)$ and $\O_{\fo} (\Xb,E^*)$ are
the space of global sections of $E^*$ of order $\l$,
respectively of finite order which is
independent of the choice of compactification $\Xb$.
Finally, the use of $\ltil$ instead of $\htil$ may seem
awkward but, in doing so, we do not deviate from our
aim of studying finite order objects. Our choice of metric
is dictated by technical reasons which will become clear in
the proof of the next theorem. So here is the main result
of this section:

\noindent \\
{\bf (5.12) Theorem:}
\emph{Let $X$ be a special affine
variety and $E$ a holomorphic vector bundle on $X$ equipped
with a Finsler metric $h$ of order $\l$. Then there is $\m \ge \l$
such that for all $q\ge 1$, we have}
\begin{align*}
H^q (\Xb, \O_\m(E^*))=0,\quad H^q (\Xb, \O_{\fo}(E^*))=0.
\end{align*}
 
\noindent
\emph{Proof:} Let $\m$ be the largest between $\l$ and the order of $g$.
For conciseness of notation we write $Y$
instead of $\P(E)$. It is known, see for instance
\cite{sommese}, that we have the isomorphism
\begin{align*}
K_Y \isom \L^{-r-1} \tensor \p^*(\text{det}(E^*))
\tensor \p^*(K_X).
\end{align*}
On the canonical line bundle $K_X$ we have the metric $k$
induced by $\f$ as follows: if, relative to a local
coordinates system $(z^1,\ldots,z^n)$, we can write
\begin{align*}
\f = \i \sum_{i,j=1}^n \f_{i\jb}\ dz^i \wedge d\zb^j
\end{align*}
then, relative to the frame
$dz^1 \wedge \ldots \wedge dz^n$ of $K_X$
we have
\begin{align*}
k= \text{det}(\f_{i\jb})^{-1}.
\end{align*}
Notice that $c_1 (K_X,k) = \text{Ric}(\f)$. Likewise, let
$\kt$ be the metric on $K_Y$ induced by $\ft$. Using the above
isomorphism we put another metric on $K_Y$ by
\begin{align*}
\kt'= \htil^{-r-1} \cdot \text{det}(G_{i\jb})^{-1} \cdot
\p^*(k).
\end{align*}
We claim that $\kt$ and $\kt'$ are almost equivalent:
there is $C >0$ such that
\begin{align*}
\tag{5.13}
C\ \r^{-n\l -nc} \ \kt' \le \kt \le
C\ \r^{-n\l} \ \kt'.
\end{align*}
We check this at a point $(x,[v])$. We may assume that the
frame of $E$ is normal at $(x,v)$. Taking into account that
det$(G_{i\jb})=1$ at $(x,[v])$ and taking determinants in
(4.11) we obtain
\begin{align*}
\text{det}(\o) \text{det}(G_{i\jb}) 2^{-n} \k^n \r^{n\l}
\text{det}(\f)
\le \text{det}(\ft) \le
\text{det}(\o) \text{det}(G_{i\jb}) 2^{n} \k^n \r^{n\l +nc}
\text{det}(\f).
\end{align*}
Recall that $\o$ is the Fubini-Study form on T$_{[v]}\P^r$.
It is known that the metric det$(\o)$ induced by $\o$ on
$\wedge^r \text{T}\P^r \isom \O_{\P^r}(r+1)$ coincides with
the canonical metric of $\O_{\P^r}(r+1)$. Therefore, under
the identification $\L_{(x,[v])}^{\tensor (r+1)}\isom
\wedge^r \text{T}_{(x,[v])}\P(E_x)$, we have $\htil^{r+1} =
\text{det}(\o)$. Taking inverse in the above inequalities
we obtain (5.13).

Under the isomorphism $\L \isom \L \tensor K_Y^{-1}
\tensor K_Y$ we put a metric $\ltil'$ on $\L$ by
\begin{align*}
\ltil'=\ltil \cdot (\kt')^{-1} \cdot \kt.
\end{align*}
By (5.13) the two metrics $\ltil$ and $\ltil'$ on $\L$
are almost equivalent:
\begin{align*}
C\ \r^{n\l} \ \ltil' \ \le \ \ltil \le C\ \r^{n\l +nc} \ \ltil'.
\end{align*}
Thus, in definition (5.11) we can replace $\ltil$ by
$\ltil'$. The latter has the advantage that its Chern form
can be bounded from below. Indeed, denoting $\L_\k$ the line
bundle $\L$ equipped with the metric
$\ltil_\k =\ltil' \cdot e^{-\k\r^\m}$,
we have:
\begin{eqnarray*}
c_1 (\L_\k)- \text{Ric}(\ft)
& = & c_1 (\L,\ltil')+ dd^c (\k \r^\m)- c_1 (K_Y, \kt) \\
& = & c_1 (\L,\ltil) - c_1 (K_Y,\kt') +dd^c (\k \r^\m) \\
& = & c_1 (\L,\htil) - c_1 (\p^*(\text{det}(E)), \text{det}(G_{i\jb}))
      + c_1 (\text{det}(E),g) \\
&   & + (r+1) c_1 (\L,\htil) + c_1 (\p^*(\text{det}(E)),
      \text{det}(G_{i\jb})) \\
&   & - \p^* \text{Ric}(\f)+dd^c(\k \r^\m) \\
& = & (r+2) c_1 (\L,\htil) + c_1 (\text{det}(E),g) \\
&   & - \p^* \text{Ric}(\f)+dd^c(\k \r^\m).
\end{eqnarray*}
From (1.2) and from the hypothesis that $g$ have order $\m$ we conclude
that there is $\k_o >0$ such that for all $\k \ge \k_o$
we have
\begin{align*}
\tag{5.14}
c_1 (\L_\k)-\text{Ric}(\ft) \ \ge \ \ft.
\end{align*}

The abstract de Rham theorem tells us that the cohomology
of a sheaf can be computed by taking an acyclic resolution.
We recall that a sheaf is said to be \emph{acyclic} if
all its cohomology groups, beside $H^0$, vanish.
Therefore, in order to show that the higher cohomology
of $\O_\m (E^*)$ vanishes, we will construct an acyclic
resolution of this sheaf which is exact at the level of
global sections. We define the sheaves $\A_\m^{0,q}(E^*)$
as follows: at each $x\in X$ the stalk $\A_\m^{0,q}(E^*)_x$
is the space of germs of smooth $\L$-valued (0,q)-forms
defined on $\p^{-1}(U)$, where $U$ is an open neighbourhood
of $x$ in $X$. If $x\in \Xb \setminus X$ is a point at
infinity, the stalk $\A_\m^{0,q}(E^*)_x$ is the space of germs
of smooth $\L$-valued (0,q)-forms defined on $\p^{-1}(U\cap X)$,
with $U$ some open neighbourhood of $x$ in $\Xb$, such that
both $u$ and $\deb u$ have order $\m$ in the L$^2$-sense.
This means that there is $\k >0$ such that
\begin{align*}
\int_{\p^{-1}(U\cap X)} |u|^2 \cdot e^{-\k \r^\m} \Phit
\ < \ \infty, \quad \quad
\int_{\p^{-1}(U\cap X)} |\deb u|^2 \cdot e^{-\k \r^\m} \Phit
\ < \ \infty.
\end{align*}
Here $|\cdot |$ is taken with respect to $\ltil'$.
Clearly $\A_\m^{0,q}(E^*),\quad 0 \le q \le n+r$, are sheaves
on $\Xb$. In fact, they are modules over the sheaf $\A_{\Xb}$
of smooth $\C$-valued functions on $\Xb$. As such they are
soft, because any $\A_{\Xb}$-module is a soft sheaf. This is
due to the fact that we can find smooth partitions of the
unity on $\Xb$. We recall that a sheaf is said to be
\emph{soft} if any section over a closed subset can be
extended to a global section. One knows that soft sheaves are
acyclic. Notice that we have a complex
\begin{align*}
\tag{5.15} 0 \lra \O_\m (E^*) \lra \A_\m^{0,0}(E^*)
\stackrel{\deb}{\lra} \A_\m^{0,1} (E^*)
\stackrel{\deb}{\lra} \ldots \stackrel{\deb}{\lra}
\A_\m^{0,q} (E^*) \stackrel{\deb}{\lra} \ldots
\end{align*}
which is clearly exact at the level of $\A_\m^{0,0}(E^*)$.
Thus the theorem will be proven once we manage to
establish that

\begin{enumerate}
\item[(i)] the complex (5.15) is exact,
\item[(ii)] the complex (5.15) is exact at the level of global
sections.
\end{enumerate}

We begin with the latter. Choose $u\in \A_\m^{0,q} (\Xb,E^*)$ such that
$\deb u=0$. There is $\k \ge \k_o$ so large that $u$ be square
integrable with respect to the metric of $\L_\k$. But (5.14)
tells us that $\L_\k$ satisfies the hypotheses of theorem (5.4).
We can find a smooth $v\in \text{L}^{0,q-1}(Y,\L_\k)$ such
that $\deb v=u$. By definition $v\in \A_\m^{0,q-1}(\Xb,E^*)$.
This proves (ii).

We now turn to (i). Exactness of a complex of sheaves means
exactness at the level of stalks. Fix an arbitrary point
$x\in \Xb$ and a germ $u_x \in \A_\m^{0,q}(E^*)_x$ represented
by some $u\in \A_\m^{0,q} (U,E^*)$. Here $U$ is a small open
neighbourhood of $x$ in $\Xb$. Our aim is to find
$v\in \A_\m^{0,q} (W,E^*)$, defined over a possibly smaller
neighbourhood $W$ of $x$, such that $\deb v=u$ on $\p^{-1}(W\cap X)$.

Let us choose a Stein neighbourhood of $x$,
say an open polycylinder $P\subset U$.
After possibly shrinking $P$ we may assume that there is
$\k \ge \k_o$ such that $u$ is square integrable relative
to $\ltil_\k$, i.e.
$u\in \text{L}^{0,q}(\p^{-1}(P\cap X), \L_\k)$.
We will fix this $\k$ for the remainder of this proof.
However, we cannot yet use (5.4) because the metric
$\ft$ on $\p^{-1}(U\cap X)$ is not complete.
We correct this by adding a horizontal term to $\ft$.

Let $\chi$ be a strictly plurisubharmonic exhaustive function
of $P$. Such a function exists because $P$ is Stein.
We put $Z=\p^{-1}(P\cap X)$. We claim that the (1,1)-form
\begin{align*}
\ft_Z = \ft + dd^c (\chi \circ \p)
\end{align*}
determines a complete metric on $Z$. To justify this we need
to show that any ball $B$ in $Z$ is relatively compact. 
Here $B$ is a ball relative to the geodesic distance induced
by $\ft_Z$. Since $\ft_Z$ dominates $\ft$ it is clear
that $B$ is included in a ball relative to the geodesic distance
induced by $\ft$. Thus $\p (B)$ stays away from the divisor at
infinity $\Xb \setminus X$. Also $\p (B)$ is included in a
ball relative to the geodesic distance induced by $dd^c\chi$
on $P$. But $dd^c \chi$ induces a complete metric on $P$
because $\chi$ is exhaustive. Thus $\p (B)$ stays away
from the boundary of $P$. We conclude that the closure of $B$
in $Z$ is compact, which justifies the claim.

Let $\kt_Z$ be the metric on $K_Z$ induced by $\ft_Z$.
It is easy to see that there are smooth functions
$\a, \b :P\lra (0,\infty)$ such that on $P\cap X$ we have
\begin{align*}
\begin{array}{rclcl}
\a \r^\l \f & \le & \frac{\k \l}{2} \r^\l \f + dd^c \chi,\\
            &     & 2 \k \r^{\l +c}\f + dd^c \chi & \le &
                             \b \r^{\l +c}\f.
\end{array}
\end{align*}
From this and (4.10) we get
\begin{align*}
\a \r^\l \f \ \le \
-\Th(v) + dd^c(\k \r^\l) + dd^c \chi
\ \le \ \b \r^{\l +c}\f.
\end{align*}
From this we obtain the analogue of (4.11):
\begin{align*}
\o + \a \r^\l \f \ \le \ \ft_Z \ \le \ \o + \b \r^{\l +c}\f.
\end{align*}
This can be used to get the analogue of (5.13):
\begin{align*}
\tag{5.16}
\b^{-n} \r^{-n\l -nc} \kt' \le \kt_Z \le
\a^{-n} \r^{-n\l} \kt'.
\end{align*}

Under the isomorphism $\L_{|Z}\isom \L_{|Z} \tensor K_Z^{-1}
\tensor K_Z$ we put a metric $\ltil_Z'$ on $\L_{|Z}$ by
\begin{align*}
\ltil_Z' = \ltil \cdot (\kt')^{-1} \cdot \kt_Z.
\end{align*}
By (5.16) the two metrics $\ltil$ and $\ltil_Z'$ on $\L_Z$
satisfy
\begin{align*}
\tag{5.17}
\a^n \r^{n\l} \ltil'_Z \ \le \ \ltil \ \le \ \b^n \r^{n\l +nc}
\ltil_Z'.
\end{align*}
The metric $\ltil_Z'$ has the advantage that we can make it
``satisfy'' the hypothesis of theorem (5.4): let
$\g :P\lra (0,\infty)$ be a strictly plurisubharmonic
exhaustive function and let $\L_\g$ be the line bundle
$\L_Z$ equipped with the metric
\begin{align*}
\ltil_\g = \ltil_Z' \cdot e^{-\k \r^\m}\cdot e^{-\g}.
\end{align*}
We have
\begin{eqnarray*}
c_1 (\L_\g) -\text{Ric}(\ft_Z)
& = & (r+2) c_1(\L,\htil) + c_1 (\text{det}(E),g) \\
&   & -\p^* \text{Ric}(\f) + dd^c (\k \r^\m) + dd^c \g \\
& \ge & \ft + dd^c \g
\end{eqnarray*}
for $\k \ge \k_o$ so that (5.14) hold. Hence for any $\g$
growing faster than $\chi$ we will have
\begin{align*}
\tag{5.18}
c_1 (\L_\g) - \text{Ric}(\ft_Z) \ \ge \ \ft_Z.
\end{align*}

Let now $Q$ be a polycylinder containing $x$ and with
$\overline{Q} \subset P$. From (5.17) and the fact that $u$
is square integrable with respect to the metrics $\ltil \cdot
e^{-\k \r^\m}$ on $\L$ and $\ft$ on $Z$ we see that $u$ is
square integrable on $\p^{-1}(Q \cap X)$ relative to the metrics 
$\ltil_Z' \cdot e^{-\k \r^\m}$ on $\L$ and $\ft_Z$ on the
manifold.

Choosing a sequence of polycylinders $\{ Q\}$ which
exhaust $P$ we prove that we can choose $\g$ growing so fast
that $u$ be square integrable on $Z$ relative to 
the metric $\ltil_\g$ on $\L$ and $\ft_Z$ on $Z$.
In other words $u$ belonds to L$^{0,q}(Z,\L_\g)$.
The hypotheses of theorem (5.4) are fulfilled: the metric
on $Z$ is complete and K\"ahler and $\L_\g$ satisfies (5.18).
We conclude that there is a smooth
$v\in \text{L}^{0,q-1}(Z,\L_\g)$ such that $\deb v=u$.

Finally, with $Q$ as above, we see that $v$ is square integrable
on $\p^{-1}(Q\cap X)$ relative to the metrics
$\ltil \cdot e^{-\k \r^\m}$ on $\L$ and $\ft$ on the manifold.
This forces $v\in \A_\m^{0,q-1} (Q,E^*)$.
The proofs of (i) and of the theorem are finished.

\noindent \\
{\bf (5.19) Corollary:} 
\emph{Let $E$ satisfy the conditions from  theorem (5.12).
Let $\F$ be a coherent algebraic sheaf on $\overline{X}$ which
is flat at all points at infinity $x\in \overline{X} \setminus X$.
Then, for $q\ge 1$,}
\begin{align*}
H^q (\overline{X}, \F \tensor_{\O_{\overline{X}}} \O_\m (E^*))=0.
\end{align*}

\noindent 
\emph{Proof:} By Hilbert's syzygy theorem $\F$ has a finite resolution
\begin{align*}
0\lra \O_{\overline{X}}^{a_n} \stackrel{\a_n}{\lra} \O_{\overline{X}}^{a_{n-1}}
\stackrel{\a_{n-1}}{\lra} \ldots 
\lra \O_{\overline{X}}^{a_1} \stackrel{\a_1}{\lra} \F \lra 0.
\end{align*}
We claim that tensoring the above with $\O_\m (E^*)$ we obtain a resolution
for the sheaf  $\F \tensor_{\O_{\overline{X}}} \O_\m (E^*)$. Indeed, the complex
\begin{align*}
\tag{5.20}
0\lra \O_{\overline{X}}^{a_n}\tensor \O_\m (E^*)
\lra \O_{\overline{X}}^{a_{n-1}}
\tensor \O_\m(E^*) \lra \ldots \lra \F \tensor \O_\m(E^*) \lra 0
\end{align*}
is exact at every point $x\in X$ by the mere fact that $\O_\m (E^*)_x \simeq
\O^r_{\overline{X},x}$ is a free, hence flat, $\O_{\overline{X},x}$-module.
When $x\in \overline{X} \setminus X$ is a point at infinity
the hypothesis that $\F_x$ be flat and standard arguments in homological
algebra ensure that ${\mathcal Ker}(\a_i)_x$ is flat for all $1\le i \le n$.
Assembling the exact sequences
\begin{align*}
0= \text{Tor}_1^{\O_{\overline{X},x}}
({\mathcal Ker}(\a_{i-1})_x ,\ \O_\m (E^*)_x) \lra
{\mathcal Ker}(\a_i)_x \tensor \O_\m (E^*)_x \lra \\
\O^{a_i}_{\overline{X},x} \tensor \O_\m (E^*)_x
\lra {\mathcal Ker}(\a_{i-1})_x \tensor \O_\m (E^*)_x \lra 0
\end{align*}
we conclude that (5.20) is exact. By the previous theorem we have
the vanishment of cohomology
\begin{align*}
H^q (\Xb, \O_{\Xb}^{a_i} \tensor \O_\m(E^*))=0,\quad q\ge 1.
\end{align*}
The corollary now follows from standard
long exact sequences in cohomology induced by (5.20).

\noindent \\
{\bf (5.21) Corollary:}
\emph{Assume that $E$ satisfies the conditions from theorem (5.12).
Then $E^*$ is spanned at every point by global sections of order $\m$.
More precisely, given $x_1, \ldots, x_N \in X$ and $e_1 \in E^*_{x_1},
\ldots, e_N \in E^*_{x_N}$ there is a global section $\s$ of $E^*$
of order $\m$ such that $\s (x_i)=e_i$ for $1 \le i \le N$.}

\noindent \\
\emph{Proof:} Let $\I \subset \O_{\overline{X}}$ be the 
ideal sheaf of $\{ x_1, \ldots, x_N \}$.
It is a coherent algebraic sheaf on $\overline{X}$
which is flat at all $x \in \overline{X} \setminus X$.
In fact, it is locally free there. By (5.19) we have
\begin{align*}
\tag{5.22}
H^1 (\overline{X}, \I \tensor_{\O_{\overline{X}}} \O_\m (E^*)) = 0.
\end{align*}
Tensoring the exact sequence
\begin{align*}
0 \lra \I \lra \O_{\overline{X}} \lra \O_{\{ x_1, \ldots, x_N \} } \lra 0
\end{align*}
with $\O_\m (E^*)$ we obtain a sequence
\begin{align*}
0 \lra \I \tensor \O_\m(E^*) \lra \O_\m(E^*) \stackrel{\a}{\lra}
\O(E^*)_{\{ x_1, \ldots, x_N \} }
\lra 0
\end{align*}
which is exact by the same discussion as in the proof of (5.19).
Its long exact sequence in cohomology together with (5.22)
show that $\a$ is surjective at the level of global sections.
This finishes the proof because $\O(E^*)_{\{ x_1, \ldots, x_N \} }$
is a skyscraper sheaf with stalks $E^*_{x_i}$ at $x_i$ and zero
outside $\{ x_1,\ldots,x_N\}$.

\noindent \\
{\bf (5.23) Proposition:}
\emph{Assume that $E$ satisfies the
conditions from theorem (5.12). Then there are finitely many
global sections $\s_1, \ldots, \s_N$ of
$E^*$ of order $\m$ which span $E^*$ at every point.
Moreover, these sections can be chosen in such a manner that
for each $x\in X$ their differentials at $x$, written
relatively to a trivialization of $E$ at $x$, span $E_x^* \tensor$ T$_x^* X$.}

\noindent \\
\emph{Proof:} We repeat here the arguments from \S 11 in
\cite{griffiths-cornalba}.
By (5.21) there exist linearly independent global
sections $\s_1, \ldots, \s_r$ of $E^*$ of order $\m$.
In this proof $r$ is the rank of $E$. Then $\s_1 \wedge \ldots
\wedge \s_r$ is a nontrivial section of det$(E^*)$.
Its zero-set $Z$ is a proper analytic
subset of $X$. We choose points $\{ x_i \}_{i \ge 1}$
on the relative interiors  of each component of $Z$.
By (11.6) in \cite{griffiths-cornalba} we can find a global section
$\s_{r+1}$ of $E^*$ of order $\m$ such that $\s_{r+1}(x_i)
\notin$ span$\{ \s_1 (x_i), \ldots, \s_{r} (x_i) \}$ for all $i$.

Let us make this explicit:
given a diverging sequence of points
$\{ x_i\}_{i\ge 1}$ in $X$ and subspaces
$F_i \subsetneqq E^*_{x_i}$, there exists
a global section $\eta$ of $E^*$ of order $\m$
such that $\eta (x_i) \notin F_i$ for all $i$.
To see this choose $e_i \in E^*_{x_i} \setminus F_i$ and global sections
$\eta_i$ of $E^*$ satisfying:
\begin{align*}
\eta_i (x_j)= 0\ \text{for}\ j<i,\quad \eta_i (x_i) = e_i,\quad 
\int_{\P(E)} |\tilde{\eta}_i |^2 \cdot  e^{-\k \r^\m} \Phit
\ < \ \infty.
\end{align*}
Such sections exist by (5.21).
The constant $\k$ can be assumed to be the same
for all $i$ because the Vanishing Theorem and its corollary (5.21)
hold with fixed large $\k$. We skip the details.
Inductively on $i$ choose nonzero
constants $c_i$ such that
\begin{align*}
\sum_{j\le i} c_j \eta_j (x_i) \notin F_i, \quad 
\int_{\P(E)} |c_i \tilde{\eta}_i |^2 \cdot
e^{-\k \r^\m} \Phit \ < \ 2^{-i}, \quad
|c_i \eta_i | < 2^{-i}\ \  \text{on} \ K_i.
\end{align*}
Here $K_i,\ i\geq 1$, is an exhaustion of $X$ by compact subsets.
The norm $|\cdot |$ is taken relative to $\ltil$.
Clearly,
\begin{align*}
\eta =\sum_{i\ge 1} c_i \eta_i
\end{align*}
is well defined and satisfies our requirements.

Repeating this procedure we construct sections $\s_{r+1},\ldots,
\s_{2r}$ of $E^*$ of order $\m$ such that $\s_1, \ldots, \s_{2r}$
span $E^*$ at all points $x_i$.
Thus $\s_1, \ldots, \s_{2r}$ span $E^*$ at all points outside
an analytic subset included in $Z$ and of strictly smaller dimension.
Repeating this procedure we can make $Z$ to be empty proving
the first part of the proposition.

What we have proven is that for all $x\in X$ the images of $\s_1, \ldots,
\s_N$ inside the vector space  $\O(E^*)_x \tensor \O_x / \mm_x$ span.

We turn now to the second part of the proposition.
We notice that there is no canonical map $E^*\lra E^* \tensor $ T$^* X$.
However, saying that $\s_1,\ldots, \s_N$ span $E^*$ at $x$ and
$d \s_1, \ldots, d \s_N$ span $E^*_x \tensor$ T$_x^* X$
is equivalent to saying that the images of $\s_1, \ldots, \s_N$ inside
$\O(E^*)_x \tensor \O_x / \mm_x^2$ span.
Indeed, we have a non-cannonical decomposition (depending
on the coordinate system)
\begin{align*}
\O_x / \mm_x^2 = \C \oplus \mm_x / \mm_x^2 = \C \oplus \text{T}_x^* X.
\end{align*}
To prove the second assertion we need only repeat the arguments from above
with $E^*_{x_i}$ replaced by $\O(E^*)_{x_i} \tensor \O_{x_i} / \mm_{x_i}^2$.\\ \\

%
%

\section{An Immersion}

The sections $\s_0,\ldots \s_N$ from (5.23) are associated
in a canonical fashion to sections $\st_0,\ldots,\st_N$ of $\L$.
The latter span $\L$ at every point, cf. (6.1). Therefore they
induce a holomorphic map
\begin{align*}
f:\P(E)\lra \P^N,\quad f(y)=(\st_0(y);\ldots;\st_N(y)).
\end{align*}
In the next section we will show that $f$ is of finite order.
Our goal in the present section is to show that the sections from
(5.23) can be chosen in such a manner that $f$ be an immersion.
A holomorphic map $f$ is said to be an \emph{immersion} if it
is one-to-one and its Jacobian has maximal rank at every
point. Equivalently, $f$ is an immersion if its image is a
complex submanifold and $f$ is a biholomorphism onto its
image. We note that the image of $f$ does not have to be
closed and, in fact, in our situation, it will not be closed.

We begin by noting at (6.1) that already $f$ is a local immersion.
However, to get the injectivity of $f$, the arguments from
(5.23) are not sufficient. We will, instead, use
arguments from the theory of Stein spaces, more precisely,
the classical way of proving that a Stein manifold of
dimension $n$ can be embedded into $\C^{2n+1}$. Our reference
will be 5.3 in \cite{hormander}. H\"ormander's lemmas can be
applied, quite generally, to a line bundle whose sections
separate points and provide holomorphic coordinates at each
point. We have checked at (5.23) that sections
of finite order of $\L$ satisfy these two properties.
One then uses the Baire category theorem to show that the
set of sections giving a ``good'' map $f$ is dense in a certain
topology. We will not be able to apply directly the category
theorem because the topology on the space of sections
of finite order of $\L$ is not complete. However, Baire's
method of proof still works, cf. (6.5).

\noindent \\
{\bf (6.1) Remark:} The sections $\st_0, \ldots, \st_N$ span $\L$
at every point, i.e. for every $y\in \P(E)$ some $\st_i (y)$ is nonzero.
Indeed, if $y= (x,[v]), \ v\in E_x \setminus \{ 0\}$, then $<\st_i, v> =
<\s_i (x), v>$. But $\s_0 (x), \ldots, \s_N (x)$ span $E_x$, hence some
$<\st_i, v>$ must be nonzero.

Thus, as already said above, there is an induced holomorphic map
\begin{align*}
f= (\st_0, \ldots, \st_N): \P(E) \lra \P^N.
\end{align*}
We claim that $f$ is a local immersion,
i.e. for all $y\in \P(E)$ the Jacobian of $f$ at $y$ has maximal
rank $n+r$. The accepted terminology is also ``$f$ is regular at $y$''.

\noindent \\
{\bf (6.2) Remark:} Let $s_0, \ldots, s_m$ be global sections of $\L$
inducing a map $f$ from $Y$ to $\P^m$. Then $f$ is regular
at a point $y\in Y$ if and only if the images of $s_0, \ldots, s_m$
under the canonical map $\L_y \lra \L_y \tensor \O_y / \mm_y^2$ generate
the latter as a vector space. \\

Let us denote by
$\Si$ the space of global sections $\s \in \G (Y,\L)$ of order $\m$ in the
L$^2$-sense:
\begin{align*}
||\s ||^2 := \int_{\P(E)} |\s |^2_{\ltil} \cdot e^{-\k \r^\m} \Phit \ < \
\infty.
\end{align*}
We choose $\k$ and $\m$ so large that the Vanishing Theorem (5.12) hold and
its corollary (5.21) apply to $\Si$. In addition, we
choose $\k$ so large that the sections $\st_0, \ldots, \st_N$ from (5.23)
belong to $\Si$. We will fix these $\k$ and $\m$ for the remainder of this
section. For the proofs of the following two lemmas we refer to
\cite{hormander} (5.3).

\noindent \\
{\bf (6.3) Lemma:} \emph{Let $K\subset Y$ be a compact subset.
Then we can find an integer $m$ and sections
$s_0, \ldots, s_m \in \Si$ inducing a map $f:Y\lra \P^m$ which regular and
one-to-one on $K$.}

\noindent \\
{\bf (6.4) Lemma:} \emph{Let $K\subset Y$ be a compact subset. Assume that some
global sections $s_0, \ldots, s_{m+1}$ of $\L$ induce a map from $Y$
to $\P^{m+1}$ which is regular and one-to-one on $K$.
Then, if $m\ge 2(n+r)+1$, we can
find $(a_0, \ldots, a_m)\in \C^{m+1}$ arbitrarily close to the origin
such that $s_0 - a_0 s_{m+1}, \ldots, s_m - a_m s_{m+1}$ induce a map
$f:Y \lra \P^m$ which is regular and one-to-one on $K$.
In fact, this is true for all $a\in \C^{m+1}$ outside a set of measure zero.}\\

We are nearing our goal. The last step is to put a topology on $\Si$ and
to show that the set of sections giving a regular one-to-one map from $Y$
to projective space is dense in $\Si^m$ equipped with the product topology.
Let $\{ K_p \}_{p\ge 1}$
be an exhaustion of $Y$ by compact subsets. We equip $\Si$ with the topology
of convergence on compact subsets: a sequence $\{ \s_q \}_{q\ge 1}$
converges to $\s$ if for each compact subset $K$ of $Y$ we have
\begin{align*}
\lim_{q\ra \infty} |\s_q - \s |_K = 0,
\end{align*}
where 
\begin{align*}
|\s_q - \s |_K =\  \text{sup}\{ |\s_q (y) - \s (y)|_{\ltil},\ y\in K \}.
\end{align*}
This topology is given by a metric invariant under translations:
\begin{align*}
d(\s', \s'') = \sum_{p\ge 1} \frac{1}{2^p} \cdot 
   \frac{|\s' - \s'' |_{K_p}}{1+|\s' - \s''|_{K_p}}.
\end{align*}
It is not clear whether $(\Si, d)$ is a complete metric space;
it is not clear that all Cauchy sequences in $\Si$ converge. However, Cauchy
sequences which are bounded in $||\cdot ||$ are convergent. Indeed,
assume that $\{ \s_q \}_{q\ge 1}$ is Cauchy relative to the distance $d$
and $||\s_q || \le M$ for some $M>0$ and all $q \ge 1$. There is a global
section $\s$ of $\L$ such that $\{ \s_q \}_{q\ge 1}$ converges pointwise
in $|\cdot |_{\ltil}$ to $\s$, uniformly on compact sets. We have
\begin{align*}
\int_{K_p} |\s_q |^2_{\ltil} \cdot e^{-\k \r^\m} \Phit \ < \ ||\s_q ||^2
                                                       \ \le \ M^2.
\end{align*}
Taking limit as $q\ra \infty$ we get
\begin{align*}
\int_{K_p} |\s |^2_{\ltil} \cdot e^{-\k \r^\m} \Phit \ \le \ M^2.
\end{align*}
Taking limit as $p\ra \infty$ we obtain $||\s || \le M$ forcing
$\s \in \Si$. Thus $\{ \s_q \}_{q\ge 1}$ converges to $\s$ in $(\Si, d)$.

\noindent \\
{\bf (6.5) Theorem:} \emph{Let $m\ge 2(n+r)+1$ be a given integer. Then the set
of (m+1)-tuples $(\s_0, \ldots, \s_m) \in \Si^{m+1}$ which induce a regular
injective map $f: Y\lra \P^m$ is dense in $\Si^{m+1}$.}

\noindent \\
\emph{Proof:} Let $G_p$ denote the set of (m+1)-tuples $(\s_0,
\ldots, \s_m) \in \Si^{m+1}$ which induce a map from $Y$ to $\P^m$ that
is regular and one-to-one on $K_p$. Clearly $G_p$ is open. We claim that
$G_p$ is dense. To see this choose an arbitrary (m+1)-tuple
$\s = (\s_0, \ldots ,\s_m) \in \Si^{m+1}$. Lemma (6.3) provides us with
sections $s_0, \ldots, s_N \in \Si$ which induce a map from $Y$ to
$\P^N$ that is regular and one-to-one on $K_p$. Then $(\s_0, \ldots,
\s_m, s_0, \ldots, s_N)$ induce a map from $Y$ to $\P^{m+N+1}$ that
is regular and one-to-one on $K_p$. Applying repeatedly lemma (6.4)
to this (m+N+1)-tuple we deduce that we can find $a_{ij},\ 0\le i \le m,
\ 0\le j\le N$ arbitrarily close to the origin such that
$\s_i' := \s_i - \sum_{j=0}^N a_{ij} s_j,\ \ 0\le i \le m$, induce a map
from $Y$ to $\P^m$ that is regular and one-to-one on $K_p$. Thus $\s'$
belongs to $G_p$ and it can be made arbitrarily close to $\s$ relative
to the distance $d$. We conclude that $G_p$ is dense.

Denoting $||\s || =$ max$\{ ||\s_j ||,\ 0\le j \le m \}$ we also notice
that $\s'$ can be made arbitrarily close to $\s$ relative to $||\cdot ||$.

We claim that $\cap_{p\ge 1} G_p$ is dense in $\Si^{m+1}$.
Had we known that $(\Si,d)$ is complete this would have been
guaranteed by Baire's category theorem. Nevertheless, in our situation,
the arguments from the proof of Baire's theorem can be carried out:

Choose arbitrary $\s \in \Si^{m+1},\ \varepsilon >0$. Choose $\s_1 \in G_1$ such
that $d(\s, \s_1 ) < \varepsilon$ and $||\s - \s_1 || < 1$.
Choose $\varepsilon_1 \in
(0, \varepsilon / 2)$ such that the closed ball $B [\s_1, 2\varepsilon_1 ]$
is contained in $G_1$.
Inductively choose $\s_p \in G_1 \cap \ldots \cap G_p$ and $\varepsilon_{p+1} \in
(0, \varepsilon_p / 2)$ such that $d(\s_p, \s_{p+1}) < 
\varepsilon_p,\ \ ||\s_{p+1}-\s_p||
< 2^{-p}$ and $B[\s_p, 2\varepsilon_p] \subset G_1 \cap
\ldots \cap G_p$. Then $\{ \s_p \}_{p\ge 1}$ is a Cauchy sequence in
$\Si^{m+1}$. Also $||\s_p || \le ||\s || +2$ for all $p$. By the discussion
preceding the theorem $\{ \s_p \}_{p\ge 1}$ converges to some $\s'
\in \Si^{m+1}$. For $q>p$ we have $d(\s_p, \s_q)\ < \ \varepsilon_p + \ldots +
\varepsilon_{q-1} \ < \ 2\varepsilon_p$. Taking limit as $q\ra \infty$ we obtain
$d(\s_p, \s') \le 2\varepsilon_p$. By construction the closed ball of center
$\s_p$ and radius $2\varepsilon_p$ is contained in $G_1 \cap \ldots \cap G_p$.
Thus $\s' \in G_1 \cap \ldots \cap G_p$. Since $p$ is arbitrary we
conclude that $\s' \in \cap_{p\ge 1} G_p$. We have $d(\s, \s') \le 2\varepsilon$.
Since $\varepsilon > 0$ is arbitrary we conclude that $\cap_{p\ge 1} G_p$
is dense in $\Si^{m+1}$. Q.e.d. \\

%
%

\section{Nevanlinna Theory on $\P(E)$}

Let $X$ be a special affine variety of dimension $n$ and $E$ a holomorphic
vector bundle of rank $r+1$ equipped with a Finsler metric of order $\l$
as defined at (5.5). Our goal is to provide means for measuring
growth of holomorphic maps from $\P(E)$ to projective space or growth
of analytic subsets of $\P(E)$. We will develop a theory along the same
lines as in section 2.

For brevity we will write $Y=\P(E)$. We also write $\r, \f,\psi$ etc.
instead of $\r \circ \p, \p^* \f, \p^* \psi$. Given $r\in (0,\infty)$
we put
\begin{align*}
Y[r]= \{ (x,[v])\in Y,\ \r'(x)\le r\},\quad
Y<r>= \{ (x,[v])\in Y,\ \r'(x)=r\}.
\end{align*}
By Sard's theorem the sets $Y<r>$ are smooth for all $r$ outside a set
of measure zero. In the sequel, each time we integrate over
$Y<r>$, it will be tacitly assumed that the latter is smooth.

We will consider $\P(E)$ equipped with the K\"ahler metric induced
by $\ft$. Our first observation is that the volume of the
balls $Y[r]$ grows at most polynomially in $r$. Indeed,
\begin{eqnarray*}
\text{vol}(Y[r])  =  \int_{Y[r]} \ft^{n+r}
                 & = & \int_{Y[r]} {n+r \choose r} (\ft^\H)^n
                          \wedge (\ft^\V)^r \\
                & \le & \int_{Y[r]} {n+r \choose r}(2\k \r^{\l +c})^n \f^n
                          \wedge (\ft^\V)^r \quad \quad
                           \text{by (5.8).}
\end{eqnarray*}
But $\f^n$ is a volume form on $X$. Thus we can apply Fubini's
theorem: the integral from above can be computed by first
integrating vertically along the fibers of $\P(E)$ and then
horizontally along $X$. We get
\begin{eqnarray*}
\text{vol}(Y[r]) & \le & \int_{x\in X[r]} {n+r \choose n}
                           (2\k \r^{\l +c})^n \f^n
                         \int_{\P(E_x)} (\ft^\V_{|\P(E_x)})^r \\
                 & =   & \int_{x\in X[r]} {n+r \choose n}
                           (2\k \r^{\l +c})^n \f^n \ \le \ r^a
\end{eqnarray*}
for some positive constant $a$. The above expression has
polynomial growth because of (1.5) and of the fact that $\f$ has
finite volume.

Given a holomorphic map $f:Y \lra \P^N$ we define its
\emph{characteristic function}
\begin{align*}
T_f(r,s)= \int_s^r \frac{dt}{t} 
\int_{Y[t]} f^* \o \wedge \psi^{n-1}\wedge \ft^r,
\end{align*}
and, for $1\le k\le n$, the \emph{higher characteristic functions}
\begin{align*}
T^{(k)}_f(r,s)= \int_s^r \frac{dt}{t} 
\int_{Y[t]} f^* \o^k \wedge \psi^{n-k}\wedge \ft^r.
\end{align*}
Here $r>s>0$ are real numbers and $\o$ is the Fubini-Study
form on $\P^N$. Given an analytic subset $Z\subset Y$ of pure
dimension $k\ge r$ we define its \emph{counting function}
\begin{align*}
N_Z(r,s)= \int_s^r \frac{dt}{t} \int_{Z[t]} \psi^{k-r}\wedge \ft^r,
\end{align*}
where $Z[t]=Z\cap Y[t]$. If the image of $Z$ under the projection
$\pp \circ \p :Y\lra \C^n$ does not contain the origin then
$N_Z (r,0)$ is well defined and we write $N_Z (r)=N_Z (r,0)$.
Given a global section $\s$ of $\O_{\P^N}(1)$ we define the
\emph{proximity function} of $f$ to the zero-set of $\s$
\begin{align*}
m_\s(r)= \int_{Y<r>} \text{log}\frac{1}{|\s \circ f|} d^c \t'
\wedge \psi^{n-1} \wedge \ft^r.
\end{align*}
Here $|\cdot |$ is taken relative to the canonical metric of
$\O_{\P^N}(1)$. Notice that $d^c \t' \wedge \psi^{n-1} \wedge \ft^r$
is a volume form on $Y<r>$. By choosing $\s$ to have norm less
than 1 at all points we can arrange that $m_\s (r)$ be non-negative.

The First Main Theorem and the Crofton Formula hold
also in this context with virtually
the same proofs as in the classical case:

\noindent \\
{\bf (7.1) Theorem:} 
\emph{Let $f:\P(E)\lra \P^N$ be a holomorphic map.
Let $H\subset \P^N$ be a hyperplane defined as the zero-set of
a global section $\s$ of $\O_{\P^N}(1)$ and $Z=f^*H$.
Assume that the image of $f$ is not contained in $H$.
Then, for $r>s>0$, we have}
\begin{align*}
N_Z(r,s)+m_\s(r)-m_\s(s)=T_f(r,s).
\end{align*}

\noindent \\
{\bf (7.2) Corollary:}
\emph{Fix $s>0$. Then, for $r>s$, and under the
hypotheses of the previous theorem, we have}
\begin{align*}
N_Z(r,s)\le T_f(r,s) +O(1).
\end{align*}
\emph{Here $O(1)$ is a constant that may depend on $s$.}

\noindent \\
{\bf (7.3) Crofton Formula:}
\emph{Let $f:\P(E)\lra \P^N$ be a holomorphic map
which is non-degenerate in the sense that the preimage of a plane $P$ of
codimension $k$ in $\P^N$ is an analytic subset of pure codimension $k$ in $\P(E)$.
Then, for $r>s>0$, we have}
\begin{align*}
T_f^{(k)} (r,s)= \int_{P\in G(N,k)} N_{f^* P} (r,s) d\m (P).
\end{align*}
\emph{Here $G(N,k)$ is the Grassmannian of planes of codimension $k$ in $\P^N$
while $\m$ is the measure defined before (2.3).}

\noindent \\
{\bf (7.4) Definition:} Let $\l$ be a non-negative real number.
Let $f: \P(E) \lra \P^N$ be a holomorphic map. We say that
\emph{$f$ has order $\l$} if there is $\k \ge 0$ such that for
some $s$ and all $r>s$ we have
\begin{align*}
T_f (r,s) \le \k r^\l.
\end{align*}
Likewise, we say that an analytic subset $Z\subset \P(E)$ of pure
dimension $k\ge r$ \emph{has order $\l$} if
\begin{align*}
N_Z (r,s) \le \k r^\l.
\end{align*}

\noindent
{\bf (7.5) Remark:} Assume that $f$ has finite order and is linearly
non-degenerate. Then (7.2) tells us that the preimage under $f$ of
any hyperplane $H\subset \P(E)$ is a divisor of finite order.
Conversely, assume that the pull-backs $f^{*}(H)$ have finite
order in a uniform manner, i.e. there are $r_o, \k, \l \ge 0$
such that for all $r \ge r_o$ and $s < r$ we have $N_{f^{*}(H)}(r,s)
\le \k r^\l$. Then, by (7.3), $f$ has finite order, too. \\

Before proceeding further let us notice that for an analytic subset
$Z\subset \P(E)$ of pure dimension $k \ge r$ we can define a counting
function
\begin{align*}
\hat{N}_Z (r,s) = \int_s^r \frac{dt}{t} \int_{Z[t]} \f^{k-r} \wedge
\ft^r
\end{align*}
by replacing $\psi$ with $\f$. One can see that $\hat{N}_Z$ and $N_Z$ have
roughly the same growth, cf. \cite{maican}.
What we mean is that there are constants
$a,b \ge 1$ independent of $Z$ such that
\begin{align*}
N_Z (r) \le r^a \hat{N}_Z (r^b),\qquad \hat{N}_Z (r) \le r^a N_Z (r^b).
\end{align*}
In particular, both counting functions have polynomial growth at the
same time. The reader may wonder why we didn't use $\hat{N}_Z$ to
measure growth in the first place. The problem with this is that
it doesn't satisfy the First Main Theorem.

\noindent \\
{\bf (7.6) Remark:} The theory we build in this section encompasses
the theory we built in section 2. What we claim is the following:
let $Z\subset X$ be a pure k-dimensional analytic subset and
$\Zt = \p^{-1}Z$. Then $Z$ and $\Zt$ have roughly the same growth.
In particular, both have finite order at the same time.
Indeed, by Wirtinger's theorem the restriction of $\f^k$ to $Z$
is a volume form on $Z$. Thus we can apply Fubini's theorem in
the definition of $\hat{N}_{\Zt}$: the integral of $\f^k \wedge
\ft^r$ on $\Zt$ can be computed by first integrating vertically
along the fibers of $\Zt$ and then horizontally along $Z$. We get
$\hat{N}_{\Zt} (r,s)= \hat{N}_Z (r,s)$. From \S 2 we know that
$\hat{N}_Z$ has roughly the same growth at $N_Z$. From the comments
preceding this remark we know that
$\hat{N}_{\Zt}$ has roughly the same growth at $N_{\Zt}$. 
This finishes the proof of the claim.

\noindent \\
{\bf (7.7) Proposition:}
\emph{Assume that $E$ is equipped with a Finsler
metric of order $\l$. Let $\s$ be a global section of $E^*$ of order
$\l$ and $Z\subset \P(E)$ the zero-set of $\st$. Then there is
$\m \ge \l$ depending only on $\l$ and $E$ such that $Z$ has order $\m$.}

\noindent \\
\emph{Proof:} Integrating the Poincar\'e-Lelong formula
\begin{align*}
[Z] = c_1 (\L,\ltil) + dd^c \text{log}|\st |^2_{\ltil}
\end{align*}
and taking logarithmic average we get
\begin{eqnarray*}
N_Z (r,s) & = & \int_s^r \frac{dt}{t}
                 \int_{Z[t]} \psi^{n-1} \wedge \ft^r \\
          & = & \int_s^r \frac{dt}{t} \int_{Y[t]}
                  \big[ c_1 (\L,\ltil)+dd^c \text{log}|\st |^2 \big]
                    \wedge \psi^{n-1} \wedge \ft^r \\
          & = & \int_s^r \frac{dt}{t} \int_{Y[t]}
                  c_1 (\L,\htil) \wedge \psi^{n-1} \wedge \ft^r \\
          &   & -\int_s^r \frac{dt}{t} \int_{Y[t]}
                  c_1 (\p^* \text{det}(E),\text{det}(G_{i\jb}))
                        \wedge \psi^{n-1} \wedge \ft^r  \\
          &   & +\int_s^r \frac{dt}{t} \int_{Y[t]}
                  c_1 (\text{det}(E),g) \wedge \psi^{n-1} \wedge \ft^r \\
          &   & +\int_s^r \frac{dt}{t} \int_{Y[t]}
                  dd^c \text{log}|\st |^2  \wedge \psi^{n-1} \wedge \ft^r.
\end{eqnarray*}
Let us denote by (i), (ii), (iii) and (iv) the integrals from the
right-hand-side above.
First we notice that (i) has polynomial growth: by (5.7) and (1.6)
\begin{eqnarray*}
(i) & \le & \int_s^r \frac{dt}{t} \int_{Y[t]}
                \ft \wedge \psi^{n-1} \wedge \ft^r \\
    & \le & \int_s^r \frac{dt}{t} \int_{Y[t]} \r^{c(n-1)} \f^{n-1}
                \wedge \ft^{r+1} \\
    & \le & \int_s^r \frac{dt}{t} \int_{Y[t]} \r^{c(n-1)} \ft^{n+r} \\
    & \le & r^{c(n-1)c'} \text{vol}(Y[r]) \ \le \ r^{a+c(n-1)c'}.
\end{eqnarray*}
Let $\m$ be the largest between $\l$ and the order of $g$.
By hypothesis $|c_1 (\text{det}(E),g)| \le \r^\m \f$ yielding
\begin{eqnarray*}
(iii) & \le & \int_s^r \frac{dt}{t} \int_{Y[t]}
                \r^\m \f \wedge \psi^{n-1} \wedge \ft^r \\
      & \le & \int_s^r \frac{dt}{t} \int_{Y[t]} \r^{\m +c(n-1)} \ft^{n+r} \\
      & \le & r^{a+\m c'+ c(n-1)c'}.
\end{eqnarray*}
By type considerations $\psi^{n-1} \wedge (\ft^\H)^2 =0$ forcing
\begin{eqnarray*}
\psi^{n-1} \wedge \ft^r & = & \psi^{n-1} \wedge (\ft^\H + \ft^\V)^r \\
& = & \psi^{n-1} \wedge (\ft^\V)^r + r\ \psi^{n-1} \wedge \ft^\H
\wedge (\ft^\V)^{r-1},\\
c_1 (G_{i\jb})\wedge \psi^{n-1} \wedge \ft^r
& = & c_1 (G_{i\jb})^\H \wedge \psi^{n-1} \wedge (\ft^\V)^r \\
&   & + r\ c_1 (G_{i\jb})^\V \wedge \psi^{n-1} \wedge \ft^\H
      \wedge (\ft^\V)^{r-1}.
\end{eqnarray*}
By (5.9), (5.10) and (5.8) we obtain
\begin{eqnarray*}
(ii) & \le & \int_s^r \frac{dt}{t} \int_{Y[t]}
                \big[ \k \r^\l \f \wedge \psi^{n-1}\wedge (\ft^\V)^r
                  + \k \r^{\l} (\ft)^\V \wedge \psi^{n-1} \wedge \ft^\H
                     \wedge (\ft^\V)^{r-1} \big] \\
     & \le & 2\k r^{a+ \l c' + c(n-1)c'}.
\end{eqnarray*}

Performing the standard ``integration twice'' procedure from Nevanlinna
theory we obtain
\begin{eqnarray*}
(iv) & = & \int_{Y<r>} \text{log}|\st | \cdot d^c \t' \wedge
              \psi^{n-1} \wedge \ft^r \\
     &   & - \int_{Y<s>} \text{log}|\st | \cdot d^c \t' \wedge
              \psi^{n-1} \wedge \ft^r.
\end{eqnarray*}
We fix $s$. By using the concavity of the logarithmic function
we obtain the estimate
\begin{align*}
\tag{7.8}
\int_{Y<r>} \text{log}|\st | \cdot d^c \t' \wedge \psi^{n-1}
\wedge \ft^r \ \le \ \text{vol}(r) \cdot \text{log}
\Big\{ \frac{1}{\text{vol(r)}} \int_{Y<r>} |\st | \cdot d^c \t' \wedge
\psi^{n-1} \wedge \ft^r \Big\}
\end{align*}
where
\begin{align*}
\text{vol}(r) = \int_{Y<r>} d^c \t' \wedge \psi^{n-1} \wedge \ft^r.
\end{align*}
But $d^c \t' \wedge \psi^{n-1}$ is a volume form on $X<r>$ away
from the branching set of the projection $\pp : X\lra \C^n$.
This set has measure zero so we ignore it. Therefore, we can apply
Fubini's theorem on $Y<r>$. We get
\begin{align*}
\tag{7.9}
\text{vol}(r) = \int_{X<r>} d^c \t' \wedge \psi^{n-1}
= \int_{\C^n <r>} d^c \text{log} ||z||^2 \wedge
(dd^c \text{log} ||z||^2)^{n-1} = 1.
\end{align*}
It remains to estimate the integral
\begin{align*}
v(r) = \int_{Y<r>} |\st | \cdot d^c \t' \wedge \psi^{n-1} \wedge \ft^r.
\end{align*}
We will first estimate the integral
\begin{align*}
w(r)=\int_1^r \frac{v(t)}{t} dt.
\end{align*}
We have
\begin{eqnarray*}
w(r) & =   & \int_{Y[1,r]} |\st |
              \cdot d\t' \wedge d^c \t' \wedge \psi^{n-1} \wedge \ft^r\\
     &     & \text{(by Fubini's theorem)} \\
     & \le & \int_{Y[1,r]}
                |\st |\cdot \r^c \f \wedge (\r^c \f)^{n-1} \wedge \ft^r\\
     &     & (\text{by (1.6) and (1.7)}) \\
     & \le & r^{c n c'} \int_{Y[1,r]} |\st | \cdot \f^n \wedge \ft^r  \\
     &     & \text{(here $c'$ is such that $\t \le c' \t'$, cf. (1.5))} \\
     & \le & r^{c n c'} \int_{Y[1,r]} |\st | \cdot \Phit \\
     & \le & r^{c n c'} 
               \Big( \int_{Y[r]} |\st |^2 \cdot e^{-\k \r^\l} \cdot \Phit \Big)^{1/2}
               \cdot \Big( \int_{Y[r]} e^{\k \r^\l} \cdot \Phit \Big)^{1/2} \\
     &     & (\text{by H\"older's Inequality}) \\
     & \le & r^{c n c'} C^{1/2} e^{\k r^\l /2} \text{vol}(Y[r])^{1/2}\\
     & \le & r^{c n c'} C^{1/2} e^{\k r^\l /2} r^{a/2}.
\end{eqnarray*}
Here $\k$ is so large that $\s$ be square integrable with respect
to the metric $\ltil \cdot e^{-\k' \r^\l}$ on $\L$ and
$\ft$ on $\P(E)$. Next we notice that
$w(r)$ is an increasing function of class $\mathcal{C}^1$ with
$w'(r)= v(r)/r$. Without loss of generality we may assume that
$\lim_{r\ra \infty}w(r)= \infty$. We claim that for all $r\ge 1$
outside a set of finite Lebesgue measure we have
\begin{align*}
w'(r) \ \le \ w(r)^2.
\end{align*}
Indeed, let $B$ be the set of ``bad'' $r$ for which $w'(r)>w(r)^2$.
Let $r_o$ be such that $w(r_o)=1$. For $r\ge r_o$ we have
\begin{align*}
\text{measure}(B\cap [1,w(r)]) < \int_1^{w(r)}
\frac{w'(t)}{w(t)^2} dt = \int_{r_o}^r \frac{dt}{t^2} < \frac{1}{r_o}.
\end{align*}
By letting $r$ converge to infinity we get that
the measure of $B$ is at most $1/r_o$
which proves the claim.

We conclude that for all $r>0$ outside a set of finite Lebesgue
measure we have
\begin{align*}
v(r) \le r^{1+ 2cnc'+a} C \cdot e^{\k r^\l}.
\end{align*}

Finally, taking logarithm and using (7.8) and (7.9) we arrive at the estimate
\begin{align*}
(iv) \le O(\text{log}(r)) + \k r^\l.
\end{align*}
We deduce the proposition holds with
$\m$ depending on the order of $g$, on $\l$ and some constants
but with $r$ outside a set of Lebesgue measure $1/r_o$.
But $T_f(r)$ is increasing hence $T_f(r) = O((r+1/r_o)^\m)=
O(r^{\m + \ep})$ with arbitrary positive $\ep$. Q.e.d.

\noindent \\
{\bf (7.10) Theorem:}
\emph{Let $X$ be a special affine variety of
dimension $n$. Let $E$ be a holomorphic vector bundle of rank $r+1$
on $X$. Assume that $E$ is equipped with a Finsler metric of order $\l$.
Then there is $\m \ge \l$ such that for any
integer $N\ge 2(n+r)+1$ there is a holomorphic immersion}
\begin{align*}
f: \P(E) \lra \P^N
\end{align*}
\emph{of order $\m$ satisfying}
\begin{align*}
f^*\O_{\P^N}(1)\isom \L.
\end{align*}

\noindent
\emph{Proof:} Let $(\s_0,\ldots, \s_N)$ be as in (6.5). Let $\m$ be
as in (7.7). For any linear combination $\s = a_0 \s_0 + \ldots +
a_N \s_N,\quad |a_i |\le 1$, the zero-set of $\st$ has order $\m$.
In fact, going through the estimates of (7.7) one sees that these
zero-sets have order $\m$ in a uniform manner. The theorem now
follows from remark (7.5). \\ \\

\newpage

%
%

\section*{List of Notations}

\begin{tabular}{rl}

$X$ & special affine variety of dimension $n$ \\
$\O_{\fo}$ & sheaf of holomorphic functions on $X$ of finite order, cf. (2.6) \\
$\tau,\ \rho=e^{\tau /2}$ & strictly plurisubharmonic exhaustive
functions on $X$, cf. (1.3) \\
$\t',\ \r'=e^{\t' /2}$ & plurisubharmonic exhaustive functions on $X$ \\
$d^c$ & the twisted derivative $\frac{1}{4\pi \i}(\de - \deb)$ \\
$dd^c$ & the complex Hessian (also known as Levi form) $\frac{1}{2\pi \i}\de \deb$ \\
$\f = dd^c \t$ & K\"ahler form on $X$ \\
$\Phi$ & the volume-form $\f^n$ of $\f$ \\
$\psi = dd^c \tau'$ & semipositive-definite form on $X$ satisfying $\psi^n=0$ \\
$\mathfrak{p}$ & generic projection from $X$ to $\mathbb{C}^n$ \\
$E$ & holomorphic vector bundle of rank $r+1$ on $X$ \\
$E^*$ & the dual of $E$ \\
$\O_{\fo}(E^*)$ & sheaf of holomorphic sections of finite order of $E^*$, cf. (5.11) \\
$\mathbb{P}(E)$ & projectivization of $E$ \\
$E_o$ & complement of the zero-section in $E$ \\
$\V$ & the vertical tangent bundle inside T$E_o$ \\
$\H$ & the horizontal tangent bundle inside T$E_o$ \\
$p$ & projection from $E$ onto the base $X$ \\
$\pi$ & projection from $\mathbb{P}(E)$ onto the base $X$ \\
$\mathcal{L}$ & the hyperplane line bundle $\mathcal{O}_{\mathbb{P}(E)}(1)$ \\
$\st$ & section of $\L$ corresponding to section $\s$ of $E^*$, cf. (4.7) \\
$g$ & metric of finite order on det$(E)$ \\
$h$ & Hermitian or Finsler metric on $E$ \\
$\tilde{h}$ & Hermitian metric on $\mathcal{L}$ induced by $h$ \\
$(G_{i\jb})$ & Hermitian metric on $\V$ induced by $h$ \\
$\ltil$ & the Hermitian metric $\htil \cdot \text{det}(G_{i\jb})^{-1}
\cdot \p^* g$ on $\L$ \\
$\ft$ & K\"ahler form on $\P(E)$, cf. (5.7) \\
$\Phit$ & the volume-form $\ft^{n+r}$ of $\ft$ \\
$\ft^\V, \ft^\H$ & the vertical, resp. the horizontal part of $\ft$,\\
& cf. paragraph after (3.23)

\end{tabular}


\begin{thebibliography}{30}

\bibitem{abate} Marco Abate and Giorgio Patrizio.
\emph{Fisler Metrics- A Global Approach},
Lecture Notes in Mathematics {\bf 1591}, Springer Verlag, Berlin, 1994.

\bibitem{andreotti-vesentini} Aldo Andreotti and Eduardo Vesentini,
\emph{Carleman Estimates for the
Laplace-Beltrami Equation on Complex Manifolds}, Publ. Math. I.H.E.S. {\bf 25}
(1965), 81-130.


\bibitem{cao-wong} Jianguo Cao and Pit-Mann Wong,
\emph{Geometry of Projectivized Vector
Bundles}, Journal of Math. of Kyoto Univ. {\bf 43}, no. 2, (2003), 369-410.


\bibitem{griffiths-carlson} James Carlson and  Phillip Griffiths,
\emph{The Order Function
for Entire Holomorphic Mappings},
Value Distribution Theory, Part A, R. O. Kujala,
A. L. Vitter (eds.), Marcel Dekker Inc., New York, 1974.

\bibitem{griffiths-cornalba} 
Maurizio Cornalba and  Phillip Griffiths, \emph{Analytic Cycles and Vector
Bundles on Non-compact Algebraic Varieties},
Invent. Math. {\bf 28} (1975), 1-106.

\bibitem{cornalba-shiffman} Maurizio Cornalba and  Bernard Shiffman,
\emph{A Counterexample
to the Transcendental Bez\^out Problem},
Ann. of Math. {\bf 96} (1972), 402-406.


\bibitem{griffiths} Phillip Griffiths, \emph{Function Theory of Finite Order on
Algebraic Varieties}, Journal of Diff. Geom. {\bf 6} (1972), 285-306 and
{\bf 7} (1972), 45-66.

\bibitem{hormander-estimates} Lars H\"ormander, \emph{L$^2$-estimates and
Existence Theorems for the $\deb$-operator}, Acta. Math. {\bf 113} (1965),
89-52.

\bibitem{hormander} Lars H\"ormander, \emph{An Introduction to Complex
Analysis in Several Variables}, North-Holland Publ. Comp.,
Amsterdam-London, 1973.

\bibitem{kodaira} Kunihiko Kodaira and  James Morrow,
\emph{Complex Manifolds}, Holt,
Rinehart and Winston, Inc., 1971.


\bibitem{lelong-gruman} Pierre Lelong and Lawrence Gruman,
\emph{Entire Functions of Several Complex Variables}, Springer Verlag,
Berlin, 1986.

\bibitem{maican} Mario Maican, \emph{Vector Bundles of Finite Order
on Affine Manifolds}, Thesis at the University of Notre Dame, July 2005.

\bibitem{wong-mulflur} James Mulflur, Albert Vitter, and Pit-Mann Wong,
\emph{Holomorphic Functions of Finite Order on Affine Varieties},
Duke Math. Journal {\bf 48}, no. 2, (1981), 389-399.


\bibitem{shabat} Boris Shabat, \emph{Distribution of Values of Holomorphic
Mappings}, Translations of Mathematical Monographs {\bf 61}, American
Mathematical Society, Providence, RI, 1985.


\bibitem{sommese} B. Shiffman, A. Sommese, \emph{Vanishing Theorems on Complex
Manifolds}, Progress in Mathematics {\bf 56}, Birkh\"auser Boston, Inc., 1985.

\bibitem{skoda} Henri Skoda,
\emph{Solution a la croissance du second probl\`eme
de Cousin dans $\C^n$}, Ann. Inst. Fourier, Grenoble {\bf 21}, no. 1, (1971),
11-23.

\bibitem{skoda_2} Henri Skoda, \emph{Sous-ensembles analytiques d'ordre fini ou
infini dans $\C^n$}, Bull. Soc. Math. France {\bf 100} (1972), 353-408.

%
%
%



\end{thebibliography}
\end{document}